\documentclass[a4paper,11pt]{article}
\usepackage[latin1]{inputenc}
\usepackage[francais,english]{babel}    
\usepackage{amssymb}
\usepackage{amsmath}
\usepackage{amsthm}
\usepackage{multicol}
\usepackage{graphicx}
\usepackage{color}
\usepackage[T1]{fontenc}
\usepackage{yfonts}
\usepackage{placeins}
\usepackage{array}
\usepackage{dsfont}
\usepackage{stmaryrd}
\usepackage{verbatim}

\newcommand{\sign}{\begin{flushright}
Thomas Haettel \\
Université Paris-Sud 11 \\
Département de Mathématiques \\
UMR 8628 CNRS \\
91405 Orsay \\
France \\
thomas.haettel@math.u-psud.fr
\end{flushright}}

\newtheorem{thm}{Théorème}[section]
\newcommand{\bthm}{\begin{thm}}
\newcommand{\ethm}{\end{thm}}

\newtheorem*{conj}{Conjecture}
\newcommand{\bconj}{\begin{conj}}
\newcommand{\econj}{\end{conj}}

\newtheorem*{question}{Question}
\newcommand{\bq}{\begin{question}}
\newcommand{\eq}{\end{question}}

\newtheorem*{thn}{Théorème}
\newcommand{\bthn}{\begin{thn}}
\newcommand{\ethn}{\end{thn}}

\newtheorem{defi}[thm]{Définition}
\newcommand{\bdf}{\begin{defi}}
\newcommand{\edf}{\end{defi}}

\newtheorem{qc}{Question de Cours}
\newcommand{\bqc}{\begin{qc}}
\newcommand{\eqc}{\end{qc}}

\newtheorem{rqc}{R\'{e}ponse}
\newcommand{\brqc}{\begin{rqc}}
\newcommand{\erqc}{\end{rqc}}

\newtheorem{exo}{Exercice}
\newcommand{\bex}{\begin{exo}}
\newcommand{\eex}{\end{exo}}

\newtheorem{sol}{Solution}
\newcommand{\bsol}{\begin{sol}}
\newcommand{\esol}{\end{sol}}

\newtheorem{pro}[thm]{Proposition}
\newcommand{\bpro}{\begin{pro}}
\newcommand{\epro}{\end{pro}}

\newtheorem{cor}[thm]{Corollaire}
\newcommand{\bcor}{\begin{cor}}
\newcommand{\ecor}{\end{cor}}

\newtheorem{lem}[thm]{Lemme}
\newcommand{\blem}{\begin{lem}}
\newcommand{\elem}{\end{lem}}

\newtheorem*{rmq}{Remarque}
\newcommand{\brq}{\begin{rmq} \upshape}
\newcommand{\erq}{\end{rmq}}

\newtheorem*{exe}{Exemple}
\newcommand{\bexe}{\begin{exe} \upshape}
\newcommand{\eexe}{\end{exe}}

\newtheorem*{pre}{Démonstration}
\newcommand{\bp}{\begin{pre} \upshape}
\newcommand{\ep}{\hfill \qed \end{pre}}
\newcommand{\epp}{\end{pre}}

\newcommand{\beq}{\begin{eqnarray*}}
\newcommand{\eeq}{\end{eqnarray*}}

\newcommand{\beqn}{\begin{equation}}
\newcommand{\eeqn}{\end{equation}}

\newcommand{\ben}{\begin{enumerate}}
\newcommand{\een}{\end{enumerate}}
\newcommand{\bit}{\begin{itemize} \renewcommand{\labelitemi}{$\bullet$} }
\newcommand{\eit}{\end{itemize}}

\newcommand{\bfg}{
\begin{figure}[H]
\begin{center}}
\newcommand{\efg}{
\end{center}
\end{figure}
\FloatBarrier}

\newcolumntype{M}[1]{>{\raggedright}m{#1}}

\newcommand{\df}{\emph}

\newcommand{\R}{\mathbb{R}}

\newcommand{\N}{\mathbb{N}}
\newcommand{\Z}{\mathbb{Z}}
\newcommand{\C}{\mathbb{C}}
\newcommand{\B}{\mathbb{B}}

\renewcommand{\P}{\mathbb{P}}

\renewcommand{\SS}{\mathbb{S}}

\newcommand{\bs}{\symbol{92}}
\newcommand{\Ch}{\mathcal{S}}

\newcommand{\ov}{\overline}

\newcommand{\sgn}{\operatorname{sgn}}
\renewcommand{\Im}{\operatorname{Im}}

\newcommand{\Ker}{\operatorname{Ker}}

\newcommand{\Diag}{\operatorname{Diag}}
\renewcommand{\dim}{\operatorname{dim}}
\newcommand{\Card}{\operatorname{Card}}

\newcommand{\Lie}{\operatorname{Lie}}

\newcommand{\Sp}{\operatorname{Sp}}

\newcommand{\bord}{\partial_\infty}

\newcommand{\Gr}{\mathcal{G}r}
\newcommand{\Plats}{\operatorname{Plats}}
\newcommand{\Cartan}{\operatorname{Cartan}}

\newcommand{\lb}{\llbracket}
\newcommand{\rb}{\rrbracket}
\newcommand{\ra}{\rightarrow}
\newcommand{\ral}[1]{\underset{#1}{\longrightarrow}}

\newcommand{\liml}{\lim\limits}

\renewcommand{\geq}{\geqslant}

\newcommand{\gothique}{\mathfrak}
\renewcommand{\gg}{\gothique{g}}
\newcommand{\kk}{\gothique{k}}
\newcommand{\pp}{\gothique{p}}
\renewcommand{\aa}{\gothique{a}}

\newcommand{\zz}{\gothique{z}}

\newcommand{\nn}{\gothique{n}}
\newcommand{\bb}{\gothique{b}}

\renewcommand{\ll}{\gothique{l}}
\newcommand{\hh}{\gothique{h}}
\renewcommand{\sl}{\gothique{sl}}

\renewcommand{\exp}{\operatorname{exp}}

\newcommand{\ad}{\operatorname{ad}}
\newcommand{\Ad}{\operatorname{Ad}}

\newcommand{\GL}{\operatorname{GL}}
\newcommand{\SL}{\operatorname{SL}}
\newcommand{\SO}{\operatorname{SO}}
\newcommand{\SU}{\operatorname{SU}}
\newcommand{\PGL}{\operatorname{PGL}}

\makeatletter
\def\Ddots{\mathinner{\mkern1mu\raise\p@
\vbox{\kern7\p@\hbox{.}}\mkern2mu
\raise4\p@\hbox{.}\mkern2mu\raise7\p@\hbox{.}\mkern1mu}}
\makeatother

\makeatletter
\def\maketitles{%
  \null
  \thispagestyle{empty}%
  \vfill
  \begin{center}\leavevmode
    \normalfont
    {\LARGE \@title\par}%
    \vskip 1.2cm
    {\large \@author\par}%
    \vskip 1.2cm
    {\large \@subtitle\par}%
    \vskip 0.8cm
    {\large \@date\par}%
  \end{center}%
  \vfill
  \null
  \cleardoublepage
  }

\def\date#1{\def\@date{#1}}
\def\author#1{\def\@author{#1}}
\def\title#1{\def\@title{#1}}
\def\subtitle#1{\def\@subtitle{#1}}
\makeatother

\title{Compactification de Chabauty\\de l'espace des sous-groupes de Cartan de $\SL_n(\R)$}
\author{Thomas Haettel}
\date{\today}

\begin{document}

\selectlanguage{francais}

\maketitle

\selectlanguage{francais}

\selectlanguage{english}
\begin{abstract} Let $G$ be a real semisimple Lie group with finite center, with a finite number of connected components and without compact factor. We are interested in the homogeneous space of Cartan subgroups of $G$, which can be also seen as the space of maximal flats of the symmetric space of $G$. We define its Chabauty compactification as the closure in the space of closed subgroups of $G$, endowed with the Chabauty topology. We show that when the real rank of $G$ is $1$, or when $G=\SL_3(\R)$ or $\SL_4(\R)$, this compactification is the set of all closed connected abelian subgroups of dimension the real rank of $G$, with real spectrum. And in the case of $\SL_3(\R)$, we study its topology more closely and we show that it is simply connected. 
\end{abstract}
\selectlanguage{francais}

\section{Introduction}

Soit $G$ un groupe de Lie réel ayant un nombre fini de composantes connexes, semi-simple, de centre fini et sans facteur compact, et soit $X$ l'espace symétrique de type non compact de $G$. L'espace $\Ch(G)$ des sous-groupes fermés de $G$ est muni d'une topologie naturelle (voir~\cite{chabauty}, et~\cite{harpe_chabauty} pour un très bon survol) qui en fait un espace compact. L'application d'isotropie de $X$ dans $\Ch(G)$ qui à un point associe son stabilisateur a été étudiée dans \cite{guivarch} et \cite{haettel_m2}. Elle identifie $X$ à l'espace des sous-groupes compacts maximaux de $G$, et l'adhérence de ce plongement dans $\Ch(G)$ est appelée la \df{compactification de Chabauty de $X$}. Guivarc'h, Ji et Taylor ont montré que cette compactification est isomorphe à la compactification de Satake-Furstenberg maximale de $X$ (voir~\cite{satake} et \cite{guivarch}). Si $X$ est l'immeuble de Bruhat-Tits d'un groupe semi-simple sur un corps local non archimédien, Y.~Guivarc'h et B.~Rémy ont défini de manière analogue la compactification de Chabauty de l'espace des sommets de $X$, vu comme espace de sous-groupes parahoriques du groupe des automorphismes de $X$, et ont montré en particulier que cette compactification est isomorphe à la compactification polyédrale de $X$ (voir~\cite{guivarch_remy}). Ce résultat a été étendu par P.E.~Caprace et J.~Lécureux à une large classe d'immeubles (voir~\cite{caprace_lecureux}). 

\bigskip

Nous nous proposons d'étudier ici une compactification naturelle de l'espace des plats maximaux de $X$, en considérant l'application qui à un plat maximal de $X$ associe l'unique sous-groupe de Cartan de $G$ le stabilisant. Ce plongement identifie l'espace des plats maximaux de $X$ à l'espace $\Cartan(G)$ des sous-groupes de Cartan de $G$, et son adhérence dans $\Ch(G)$, notée $\ov{\Cartan(G)}^\Ch$, est appelée la \df{compactification de Chabauty} de $\Cartan(G)$. Puisque le stabilisateur d'un plat maximal dans $G$ est le normalisateur $N_G(A)$ d'un sous-groupe de Cartan $A$ de $G$, nous étudions donc une compactification naturelle de l'espace homogène $G/N_G(A)$.

\bigskip

Lorsque $G$ est déployé sur $\R$, le sous-groupe de Cartan stabilisant un plat maximal de $X$ est la composante neutre de son stabilisateur dans $G$. De plus le sous-groupe $N_G(A)$ est également le normalisateur du tore déployé sur $\R$ maximal $Z_G(A)$, ainsi l'espace $G/N_G(A)$ s'identifie aussi à l'espace des tores déployés sur $\R$ maximaux de $G$. Par ailleurs, l'espace $G/N_G(A)$ est un quotient fini de l'espace homogène $G/A$, dont la géométrie et la dynamique restent à comprendre finement. Par exemple, une conjecture de Margulis (voir~\cite{margulis_conjecture} et \cite{maucourant}) affirme que les $A$-orbites bornées dans $\SL_3(\R)/\SL_3(\Z)$ sont compactes (où $A$ est le sous-groupe diagonal connexe). D'un point de vue dual, ceci motive l'étude de compactifications naturelles des espaces homogènes $G/A$.

\bigskip

Cette compactification de Chabauty a également été définie et étudiée pour $\SL_n(\C)$ par Iliev et Manivel (voir~\cite{iliev_manivel} et \cite{iliev_manivel_severi}), puis par Le Barbier Grünewald pour des groupes complexes plus généraux (voir~\cite{lebarbier} et \cite{lebarbier_ex}), dans la perspective d'obtenir des exemples naturels de variétés algébriques complexes aux propriétés intéressantes. Cependant, la plupart de leurs résultats est propre au cadre complexe, et ne peut être transposée directement au cadre réel. Dans le cas de $G=\SL_n(\C)$, l'espace homogène $G/N_G(A)$ s'identifie aux $n$-uplets de points génériques dans $\P^{n-1}(\C)$, dont d'autres compactifications ont été étudiées du point de vue de la géométrie algébrique complexe, notamment le schéma de Hilbert des $n$ points de $\P^{n-1}$ et la compactification de Fulton-Macpherson de l'espace des configurations de points ordonnés sur une variété (voir~\cite{fulton_macpherson}).

\bigskip

Après avoir défini formellement la compactification de Chabauty de l'espace des sous-groupes de Cartan de $G$, nous établissons plusieurs propriétés générales de $\ov{\Cartan(G)}^\Ch$. Nous montrons qu'il suffit d'étudier les sous-espaces de Cartan de l'algèbre de Lie de $G$. Nous montrons que, contrairement au cas complexe, cette compactification n'est pas Zariski-fermée dans la grassmannienne de l'algèbre de Lie de $G$, mais est seulement une sous-variété semi-algébrique réelle.

La compactification de Chabauty $\ov{\Cartan(G)}^\Ch$ est naturellement incluse dans le sous-espace ${\cal A}(G)$ de $\Ch(G)$ constitué des sous-groupes fermés abéliens connexes de dimension le rang réel de $G$ inclus dans un sous-groupe de Borel de $G$ (ou, de manière équivalente, dont le spectre de l'action adjointe est réel).

\bthn Supposons que $G$ soit un produit presque direct de groupes de rang réel un et de copies de $\SL_3(\R)$ et de $\SL_4(\R)$. Alors $\ov{\Cartan(G)}^\Ch = {\cal A}(G)$, et est la réunion d'un nombre fini d'orbites de $G$. (Voir les théorèmes~\ref{thm:sgabeliens_rg1}, \ref{thm:sgabeliens_sl3} et \ref{thm:sgabeliens_sl4}.) \ethn

Par contre, si $n \geq 7$ nous savons que l'inclusion de $\ov{\Cartan(\SL_n(\R))}^\Ch$ dans ${\cal A}(\SL_n(\R))$ est stricte (voir le lemme~\ref{lem:sgabeliens_sl7}).

Dans le cas particulier où $G$ est de rang réel un, la topologie de $\ov{\Cartan(G)}^\Ch$ est explicite.

\bthn Supposons que le rang réel de $G$ soit $1$. Notons $\bord X$ la sphère à l'infini de $X$ et $\bord X^{(2)}$ l'espace des paires de points éventuellement confondus de $\bord X$. Alors $\ov{\Cartan(G)}^\Ch$ est naturellement $G$-homéomorphe à l'éclatement de $\bord X^{(2)}$ le long de la diagonale. (Voir le théorème~\ref{thm:sgabeliens_rg1}.) \ethn

Dans le cas où $G=\SL_3(\R)$, bien que le groupe fondamental de l'espace $\Cartan(G)$ soit de cardinal $48$ (voir le lemme~\ref{lem:pi1_plats}), nous montrons le résultat suivant.

\bthn L'espace $\ov{\Cartan(\SL_3(\R))}^\Ch$ est simplement connexe. (Voir le théorème~\ref{thm:chabautysl3}.) \ethn

Je tiens à remercier très chaleureusement Frédéric Paulin, pour ses très nombreux conseils avisés et relectures attentives. Je souhaite également remercier Olivier Benoist pour de nombreuses discussions, ainsi que Laurent Manivel.

\section{Définition de la compactification}

Soit $G$ un groupe de Lie réel ayant un nombre fini de composantes connexes, semi-simple, de centre fini et sans facteur compact.

L'espace $\Ch(G)$ des sous-groupes fermés de $G$ est muni de la topologie de Chabauty (voir~\cite{chabauty} et \cite{harpe_chabauty}) qui en fait un espace compact. Puisque $G$ est métrisable, l'espace $\Ch(G)$ l'est également, et on peut décrire aisément la convergence des suites. Une suite de sous-groupes fermés $(H_n)_{n \in \N}$ de $G$ converge vers un sous-groupe fermé $H$ si et seulement si les deux conditions suivantes sont vérifiées.
\ben
\item Pour tout $h \in H$, il existe une suite $(h_n)_{n \in \N}$ convergeant vers $h$ telle que, pour tout $n \in \N$, nous ayons $h_n \in H_n$ ;
\item pour toute partie infinie $P \subset \N$, et pour toute suite $(h_n)_{n \in P}$ convergeant vers $h$ telle que $h_n \in H_n$ pour tout $n \in P$, nous avons $h \in H$.
\een

Dans cet article, un \df{sous-groupe de Cartan} de $G$ est la composante neutre $A$ d'un tore déployé sur $\R$ maximal de $G$, autrement dit un sous-groupe fermé connexe de $G$ dont l'algèbre de Lie est un \df{sous-espace de Cartan} de $\gg$, c'est-à-dire une sous-algèbre de Lie abélienne diagonalisable sur $\R$ maximale de $\gg$. Un sous-espace de Cartan de $\gg$ n'est une sous-algèbre de Cartan de $\gg$, c'est-à-dire dont la complexifiée est une sous-algèbre de Cartan de la complexifiée de $\gg$, que si $G$ est déployé sur $\R$.

Notons $\Cartan(G)$ l'espace des sous-groupes de Cartan de $G$, muni de la topologie induite par la topologie de Chabauty sur l'espace $\Ch(G)$ des sous-groupes fermés de $G$. Le groupe $G$ agit transitivement sur $\Cartan(G)$ par conjugaison. On appelle \df{compactification de Chabauty} de $\Cartan(G)$, et on note $\ov{\Cartan(G)}^\Ch$, l'adhérence de $\Cartan(G)$ dans l'espace compact $\Ch(G)$.

Fixons un sous-groupe de Cartan $A$ de $G$. Notons $\gg$ l'algèbre de Lie de $G$, $\aa$ l'algèbre de Lie de $A$, $r = \dim \aa$ le rang réel de $G$ et $\Cartan(\gg)$ l'espace des sous-espaces de Cartan de $\gg$, muni de la topologie induite par celle de la grassmannienne $\Gr_r(\gg)$ des $r$-plans de $\gg$. Le groupe $G$ agit sur $\Gr_r(\gg)$ et donc sur $\Cartan(\gg)$ par l'action adjointe. L'application qui à un sous-groupe de Cartan de $G$ associe son algèbre de Lie définit un homéomorphisme $G$-équivariant de $\Cartan(G)$ sur $\Cartan(\gg)$. On appelle \df{compactification de Chabauty} de $\Cartan(\gg)$, et on note $\ov{\Cartan(\gg)}^\Ch$, l'adhérence de $\Cartan(\gg)$ dans la grassmannienne $\Gr_r(\gg)$. Les deux compactifications ainsi définies sont isomorphes, au sens de la proposition suivante.

\bpro L'application
\beq \eta : \ov{\Cartan(G)}^\Ch & \ra & \ov{\Cartan(\gg)}^\Ch \\
H & \mapsto & \Lie(H)\eeq
est un homéomorphisme $G$-équivariant. \epro

\bp Notons $\Sp(\ell)$ le spectre d'un endomorphisme $\ell$ de $\gg$. L'application exponentielle de $G$ est un difféomorphisme de $\mathcal{E}=\{X \in \gg \,:\, \Sp(\ad X) \subset \R\}$ sur $\{g \in G \,:\, \Sp(\Ad g) \subset \R\}$. Ainsi $\eta$ est l'application qui à $H \in \ov{\Cartan(G)}^\Ch$ associe $(\exp|_\mathcal{E})^{-1}(H)$, c'est un homéomorphisme. \ep

Par ailleurs, la topologie de $\Cartan(G)$ s'identifie à la topologie de l'espace homogène $G/N_G(A)$.

\bpro \label{pro:sg_ssalg} Les applications
\beq \phi : G/N_G(A) \ra \Ch(G) & \mbox{ et } & \varphi : G/N_G(A) \ra \Gr_r(\gg) \\
gN_G(A) \mapsto gAg^{-1} & & gN_G(A) \mapsto \Ad g(\aa)\eeq
sont des plongements $G$-équivariants. \epro

\bp L'équivariance de ces deux applications est immédiate. Et puisque $\phi = \eta^{-1} \circ \varphi$, il suffit de montrer que $\varphi$ est un plongement.

Si $\Ad G$ est un sous-groupe algébrique de $\GL(\gg)$, l'application $\varphi$ s'identifie à l'application orbitale pour l'action algébrique du groupe algébrique $\Ad G$ sur la variété algébrique projective $\Gr_r(\gg)$, car le fixateur de $\aa$ pour l'action adjointe de $G$ est $N_G(A)$. C'est une immersion de variétés algébriques, or son image est localement fermée pour la topologie de Zariski (voir par exemple \cite[Proposition~8.3, p.~60]{humphreys}), c'est donc en particulier un plongement pour la topologie analytique. Si $\Ad G$ n'est pas un sous-groupe algébrique de $\GL(\gg)$, l'application $\varphi$ s'identifie à la restriction à $G/N_G(A)$ de l'application orbitale pour l'action algébrique de l'adhérence de Zariski $\ov{\Ad G}^Z$ sur la variété algébrique projective $\Gr_d(\gg)$, c'est donc également un plongement.
\ep

Soit $X$ l'espace symétrique de type non compact de $G$, notons $\Plats(X)$ l'espace des plats maximaux de $X$, muni de la topologie induite par la topologie de Chabauty sur l'espace des fermés de $X$. Soit $P_0$ l'unique plat maximal de $X$ stabilisé par $A$. Considérons le plongement
\beq \Plats(X) & \ra & \Ch(G) \\
P = g \cdot P_0 & \mapsto & gAg^{-1} \mbox{ l'unique sous-groupe de Cartan de $G$ qui stabilise $P$ },\eeq
c'est un homéomorphisme $G$-équivariant de $\Plats(X)$ sur $\Cartan(G)$, dont l'adhérence de l'image est la compactification de Chabauty de l'espace des plats maximaux de $X$. Remarquons que si $G$ est déployé sur $\R$, alors le sous-groupe de Cartan de $G$ qui stabilise $P$ est la composante neutre du stabilisateur de $P$, ce qui simplifie l'expression de cette application d'isotropie.

\section{Propriétés générales de la compactification}

Considérons l'algèbre de Lie complexifiée $\gg_\C = \gg \otimes \C$ de $\gg$. L'espace $\ov{\Cartan(\gg)}^\Ch$ se plonge dans la grassmannienne $\Gr_r(\gg_\C)$ des $r$-plans complexes de $\gg_\C$, en associant à une sous-algèbre de Lie réelle de $\gg$ son produit tensoriel avec $\C$. Notons $\ov{\Cartan(\gg_\C)}^{Zar}$ l'adhérence de Zariski complexe de ce plongement, et décrivons le lien avec la compactification de Chabauty réelle $\ov{\Cartan(\gg)}^\Ch$.

Dans le cas où $G$ est déployé sur $\R$, cette adhérence $\ov{\Cartan(\gg_\C)}^{Zar}$ est la compactification de Chabauty de l'espace des sous-algèbres de Cartan complexes de $\gg_\C$, appelée \df{variété des réductions} par Iliev et Manivel (voir~\cite{iliev_manivel} et \cite{iliev_manivel_severi}).

Remarquons que la variété algébrique complexe $\ov{\Cartan(\gg_\C)}^{Zar}$ est irréductible. En effet, c'est l'adhérence de Zariski de l'orbite de $\aa_\C = \aa \otimes \C$ par le groupe algébrique connexe $(\Ad G)_\C$, complexifié du groupe $\Ad G \subset \GL(\gg_\C)$. Le lemme suivant est immédiat.

\blem Soit $\ll \in \ov{\Cartan(\gg_\C)}^{Zar}$. Si $\ll$ est la tensorisée avec $\C$ d'un élément de $\ov{\Cartan(\gg)}^\Ch$, alors $\ll$ est définie sur $\R$ et $\Sp(\ad \ll(\R)) \subset \R$. \elem

La compactification de Chabauty $\ov{\Cartan(\gg)}^\Ch$ est l'adhérence pour la topologie analytique réelle de $\Cartan(\gg)$, qui est une sous-variété algébrique réelle de $\Gr_r(\gg)$, ainsi $\ov{\Cartan(\gg)}^\Ch$ est une variété semi-algébrique réelle (voir~\cite[Proposition~2.2.2,p.~27]{bochnak_coste_roy}). Par contre, contrairement au cas complexe, nous avons le résultat suivant.

\bpro La compactification de Chabauty $\ov{\Cartan(\gg)}^\Ch$ n'est pas Zariski-fermée dans $\Gr_r(\gg)$.\epro

\bp
Remarquons tout d'abord que ceci est vrai pour $\gg=\sl_2(\R)$. L'espace $\Cartan(\gg)$ des sous-espaces de Cartan de $\gg$ est le sous-espace de $\P(\gg)$ constitué des classes d'homothétie d'éléments de $\gg$ dont les deux valeurs propres (opposées) sont réelles et différentes de $0$.

L'adhérence $\ov{\Cartan(\gg)}^\Ch$ de $\Cartan(\gg)$ pour la topologie analytique réelle est le sous-espace de $\P(\gg)$ constitué des classes d'homothétie d'éléments de $\gg$ dont les deux valeurs propres sont réelles. En revanche, l'adhérence de $\Cartan(\gg)$ pour la topologie de Zariski réelle est tout $\P(\gg)$.

\bigskip

Dans le cas général d'une algèbre de Lie semi-simple $\gg$ de type non compact, considérons une sous-algèbre de Lie de $\gg$ isomorphe à $\sl_2(\R)$. Le même argument montre que les adhérences de Zariski réelle et analytique réelle de $\Cartan(\gg)$ sont distinctes.
\ep

Si $\Ad G$ est le groupe de Lie ajdoint de $G$, considérons le morphisme surjectif de groupes de Lie $\Ad : G \ra \Ad G$ de noyau le centre de $G$ fini. Alors la composition avec $\Ad$ définit un homéomorphisme $\Ad$-équivariant de $\ov{\Cartan(G)}^\Ch$ sur $\ov{\Cartan(\Ad G)}^\Ch$. Ainsi l'étude de la compactification de Chabauty $\ov{\Cartan(G)}^\Ch$ de l'espace $\Cartan(G)$ ne dépend que de l'algèbre de Lie $\gg$ de $G$.

\bigskip

Si l'algèbre de Lie $\gg$ est le produit direct $\gg = \bigoplus_{i=1}^n \gg_i$ de $n$ sous-algèbres de Lie semi-simples, alors l'application
\beq \prod_{i=1}^n \Cartan(\gg_i) & \ra & \Cartan(\gg) \\
(\aa_i)_{i \in \lb 1,n\rb} & \mapsto & \bigoplus_{i=1}^n \aa_i\eeq
est un homéomorphisme qui s'étend en un homéomorphisme $G$-équivariant du produit $\prod_{i=1}^n \ov{\Cartan(G_i)}^\Ch$ sur $\ov{\Cartan(G)}^\Ch$. Ainsi on peut se ramener au cas où $\gg$ est une algèbre de Lie réelle simple de type non compact.

Le premier théorème de l'introduction découle donc des parties~\ref{sec:rang1}, \ref{sec:sl3} et \ref{sec:sl4}.

\bigskip

Soit $G = KAN = KNA$ une décomposition d'Iwasawa de $G$, où $K$ est un sous-groupe compact maximal de $G$ et où $N$ est un sous-groupe unipotent de $G$. Soit $x_0 \in P_0$ le point de $X$ dont le stabilisateur soit $K$. Pour décrire les limites de sous-espaces de Cartan de $\gg$, on peut donc modulo l'action du groupe compact $K$ se ramener à étudier les limites d'images de $\aa$ par des éléments de $N$.

Il est immédiat que la compactification de Chabauty $\ov{\Cartan(G)}^\Ch$ est incluse dans le sous-espace fermé de $\Ch(G)$ constitué des sous-groupes fermés abéliens connexes de dimension $r$ de $G$.

Dans cet article, on appelle \df{sous-groupe de Borel} de $G$ le normalisateur $B$ du sous-groupe $AN$ d'une décomposition d'Iwasawa $G=KAN$ de $G$. C'est un sous-groupe algébrique résoluble, dont la composante neutre est $AN$. Son algèbre de Lie $\Lie(B) = \Lie(AN)$ est appelée une \df{sous-algèbre de Borel} de $\gg$ .

\bpro Soit $H$ un sous-groupe fermé abélien connexe de dimension $r$ de $G$. Les assertions suivantes sont équivalentes, et sont vérifiées si $H \in \ov{\Cartan(G)}^\Ch$.
\ben
\item Le sous-groupe $H$ est inclus dans un sous-groupe de Borel de $G$.
\item Le spectre de l'action adjointe de $H$ sur $\gg$ est réel.
\een
\epro

\bp Si $H$ est inclus dans un sous-groupe de Borel $B$, puisque le spectre de l'action adjointe de $B$ sur $\gg$ est réel, le spectre de $H$ l'est également. 

Réciproquement, supposons que le spectre de l'action adjointe de $H$ soit réel. Alors d'après~\cite[Theorem~3.2, p.~133]{moore_amenable}, les sous-groupes moyennables fermés connexes maximaux $L$ de $G$ forment un nombre fini de classes de conjugaison, et chacun d'eux est une extension d'un sous-groupe d'un groupe de Borel par un sous-groupe compact $K_L$. Soit $L$ un sous-groupe moyennable connexe maximal de $G$ contenant $H$, tel que le sous-groupe compact $K_L$ de $L$ soit minimal. Puisque le spectre de $H$ est réel, le sous-groupe $K_L$ de $L$ est trivial, donc $L$ est inclus dans un sous-groupe de Borel de $G$, et $H$ également.
\ep

Ainsi la compactification de Chabauty $\ov{\Cartan(G)}^\Ch$ est incluse dans le sous-espace fermé de $\Ch(G)$
\beq {\cal A}(G) &=& \{ \mbox{sous-groupes fermés de $G$ abéliens connexes de dimension r}\\
& & \mbox{vérifiant ces deux propriétés équivalentes}\}.\eeq
Notons de même ${\cal A}(\gg)$ le sous-espace de $\Gr_r(\gg)$ des sous-algèbres de Lie abéliennes incluses dans une sous-algèbre de Borel de $\gg$ ou, de manière équivalente, dont le spectre de l'action adjointe est réel.

\bigskip

La question naturelle que l'on se pose est alors de savoir dans quel cas nous avons $\ov{\Cartan(G)}^\Ch = {\cal A}(G)$. Nous montrons dans cet article que c'est le cas si $G$ est de rang réel un (voir le théorème~\ref{thm:sgabeliens_rg1}), ainsi que pour $\SL_3(\R)$ (voir le théorème~\ref{thm:sgabeliens_sl3}) et $\SL_4(\R)$ (voir le théorème~\ref{thm:sgabeliens_sl4}). Par contre, comme le font remarquer Iliev et Manivel (dans le cas complexe), ce n'est pas le cas pour $\SL_m(\R)$ si $m \geq 7$, comme l'explicite le lemme suivant.

\blem \label{lem:sgabeliens_sl7} Pour $m \geq 7$, il existe une sous-algèbre abélienne de $\sl_m(\R)$ de dimension $m-1$ incluse dans une sous-algèbre de Borel, qui ne soit pas limite de sous-espaces de Cartan. \elem

\bp Nous reprenons ici la preuve de~\cite[p.~3]{iliev_manivel} dans le cas réel. Fixons $V$ un sous-espace vectoriel réel de dimension $p=\lfloor \frac{m}{2} \rfloor$ de $\R^m$. L'ensemble $X_V$ des endomorphismes $f$ de $\R^m$ tels que $\Im f \subset V \subset \Ker f$ est un sous-espace vectoriel de dimension $p(m-p)$, c'est une sous-algèbre de Lie abélienne de $\sl_m(\R)$ incluse dans la sous-algèbre de Borel standard. Tout sous-espace vectoriel de dimension $m-1$ de $X_V$ est donc une sous-algèbre de Lie abélienne de $\sl_m(\R)$ de dimension $m-1$ incluse dans une sous-algèbre de Borel, ainsi c'est un élément de ${\cal A}(\sl_m(\R))$. Par ailleurs un sous-espace vectoriel de $X_V$ de dimension $m-1$ générique détermine uniquement $V$. L'ensemble des sous-espaces vectoriels de dimension $m-1$ de $X_V$, où $V$ parcourt les sous-espaces vectoriels de dimension $p$ de $\R^m$, est donc une sous-variété incluse dans ${\cal A}(\sl_m(\R))$ de dimension $\dim \Gr_{m-1}(p(m-p)) + \dim \Gr_p(m) = (m-1)(p(m-p)-m+1)+p(m-p)$. Cette dimension est strictement supérieure à la dimension $m(m-1)$ de $\Cartan(\sl_m(\R))$ dès que $m \geq 7$. D'après~\cite[Proposition~2.8.13]{bochnak_coste_roy}, cette sous-variété ne peut être entièrement incluse dans l'adhérence pour la topologie analytique de la sous-variété semi-algébrique $\Cartan(\sl_m(\R))$. \ep

\bigskip

Remarquons que l'espace $\Cartan(G)$ n'est pas simplement connexe, comme le précise le lemme suivant.

\blem \label{lem:pi1_plats} L'espace $\Plats(X) \simeq G/N_G(A)$ se rétracte par déformation forte sur $\{P \in \Plats(X) \,:\, x_0 \in P\} \simeq K / N_K(A)$. Ainsi son groupe fondamental est isomorphe à $W \ltimes \pi_1(K/Z_K(A))$, où $W=N_K(A)/Z_K(A)$ désigne le groupe de Weyl de $G$. En particulier si $G=\SL_n(\R)$ avec $n \geq 3$, alors $\pi_1(\Cartan(G))$ est de cardinal $n!2^n$. \elem

\bp Notons $\theta$ l'involution de Cartan de $\gg$ associée au choix du point base $x_0$, et $\pp$ le sous-espace propre de $\theta$ pour la valeur propre $-1$.

Fixons $t \in [0,1]$, soit $P$ un plat maximal de $X$ et soit $x$ le projeté orthogonal de $x_0$ sur $P$. D'après la décomposition polaire en $x_0$, soit $Y$ l'unique élément de $\pp$ tel que $\exp(Y) \cdot x_0 = x$. Posons $f_t(P) = \exp(-tY) \cdot P$.

Nous avons ainsi défini rétraction par déformation forte $f : \Plats(X) \times [0,1] \ra \Plats(X)$ de $\Plats(X)$ sur $\{P \in \Plats(X) \,:\, x_0 \in P\} \simeq K/N_K(A)$.

\bigskip

Si $G=\SL_n(\R)$ avec $n \geq 3$, alors $W \simeq \frak{S}_n$ donc $W$ est de cardinal $n!$. De plus $K \simeq \SO(n)$ et $Z_K(A)$ est isomorphe au sous-groupe diagonal de coefficients $\pm 1$. Ainsi $\Card \pi_1(K/Z_K(A)) = \Card(Z_K(A)) \times \Card \pi_1(\SO(n)) = 2^{n-1} \times 2 = 2^n$.\ep

\bigskip

Énonçons les conséquences immédiates sur $\ov{\Cartan(G)}^\Ch$ de certains résultats de A.~Iliev, L.~Manivel et M.~Le Barbier Grünewald sur $\ov{\Cartan(\gg_\C)}^{Zar}$. Rappelons qu'un élément de $\gg$ est dit \df{régulier} si son centralisateur dans $\gg$ est de dimension minimale. Une sous-algèbre abélienne $\ll$ de $\gg$ est dite \df{régulière} s'il existe un élément de $\ll$ régulier dont le centralisateur soit égal à $\ll$. Notons $\ov{\Cartan(\gg)}^{Zar}$ l'adhérence de $\Cartan(G)$ dans $\Gr_r(\gg)$ pour la topologie de Zariski (réelle), elle contient $\ov{\Cartan(G)}^\Ch$.

\bthm
\ben
\item Toute sous-algèbre de Lie régulière de $\gg$ appartient au lieu lisse de $\ov{\Cartan(\gg)}^\Ch$.
\item Toute sous-algèbre de Lie de $\gg$ appartenant à $\ov{\Cartan(\gg)}^{Zar}$ est l'algèbre de Lie d'un sous-groupe algébrique de $G$.
\item Si la $G$-orbite d'une sous-algèbre de Lie $\ll$ de $\gg$ appartenant à $\ov{\Cartan(\gg)}^{Zar}$ est fermée, alors $\ll$ est constituée d'éléments nilpotents.
\item Pour $G=\SL_3(\R)$, alors $\ov{\Cartan(\sl_3(\R))}^{Zar}$ est lisse.
\item Pour $G=\SL_4(\R)$, alors le lieu singulier de $\ov{\Cartan(\sl_4(\R))}^{Zar}$ est la réunion de deux orbites isomorphes à $\P(\R^4)$.
\een
\ethm

\bp Les démonstrations sont les m\^{e}mes, ou en découlent immédiatement, que~\cite[Theorem~3.7]{lebarbier}, \cite[Corollary~5.3]{lebarbier}, \cite[Proposition~5.5]{lebarbier} et \cite[Proposition~9]{iliev_manivel}.
\ep

Par contre, la compactification de Chabauty $\ov{\Cartan(\sl_3(\R))}^\Ch$ est un fermé strict de $\ov{\Cartan(\sl_3(\R))}^{Zar}$ donc n'est pas lisse.

\section{Le cas de rang un}

\label{sec:rang1}

Soit $G$ un groupe de Lie réel connexe semi-simple de centre fini, sans facteur compact et de rang réel un, et soit $X$ son espace symétrique (ainsi $X$ est l'espace hyperbolique réel, complexe, quaternionique ou le plan hyperbolique octonionique, voir~\cite[Theorem~9.1.1, p.~83]{parker_hyperbolic}). Notons $\bord X$ la sphère à l'infini de $X$, notons $\bord X^{(2)}$ l'espace des paires de points éventuellement confondus, et notons $\Diag(\bord X^{(2)}) \subset \bord X^{(2)}$ la diagonale.

Remarquons que l'application qui à un sous-groupe de Cartan de $G$ associe les extrémités de l'unique géodésique de $X$ translatée par le sous-groupe de Cartan est un homéomorphisme $G$-équivariant de $\Cartan(G)$ sur $\bord X^{(2)} \bs \Diag(\bord X^{(2)})$.

\bthm \label{thm:sgabeliens_rg1} La compactification de Chabauty $\ov{\Cartan(G)}^\Ch$ de l'espace des sous-groupes de Cartan de $G$ coïncide avec ${\cal A}(G)$. Sous l'action de $G$, elle est la réunion de l'orbite ouverte des sous-groupes de Cartan et de
\bit
\item l'orbite fermée des sous-groupes fermés connexes unipotents de dimension $1$, si $X$ est hyperbolique réel ;
\item deux orbites (dont une seule est fermée) de sous-groupes fermés connexes unipotents de dimension $1$, si $X$ est hyperbolique complexe, quaternionique ou le plan hyperbolique octonionique.
\eit
De plus, $\ov{\Cartan(G)}^\Ch$ est $G$-isomorphe à l'éclaté de $\bord X^{(2)}$ le long de $\Diag(\bord X^{(2)})$. \ethm

\bp  Soit $\hh$ une sous-algèbre de Lie de $\gg$ de dimension $1$ incluse dans une sous-algèbre de Borel $\bb$ de $\gg$, qui n'est pas un sous-espace de Cartan de $\gg$. \'{E}crivons $\nn = \gg_\alpha \oplus \gg_{2\alpha}$ la décomposition en espaces de racines de $\nn$ sous l'action adjointe de $A$ ($\gg_{2\alpha}=\{0\}$ si $X$ est hyperbolique réel). Soit $Y = Y_\alpha + Y_{2\alpha} \in \nn = \gg_\alpha \oplus \gg_{2\alpha}$ non nul tel que $\hh = \R Y$. Alors la suite de sous-espaces de Cartan $(\Ad \exp(2nY_\alpha +nY_{2\alpha})(\aa))_{n \in \N}$ converge vers $\hh$ dans $\ov{\Cartan(\gg)}^\Ch$. Ainsi $\ov{\Cartan(\gg)}^\Ch = {\cal A}(\gg)$.

\bigskip

D'après~\cite[Theorem~3.4.1, p.~72]{chen_greenberg}, il existe une seule classe de conjugaison d'éléments unipotents (différents de l'identité) si $G$ localement isomorphe à $\SO_0(n,1)$ avec $n \geq 2$, et deux classes si $G$ est localement isomorphe à $\SU(n,1)$ ou à $\Sp(n,1)$, avec $n \geq 2$. D'après~\cite[Theorem~4.4, p.~11]{allcock_cayley}, il existe deux classes de conjugaison d'éléments unipotents (différents de l'identité) si $G$ est localement isomorphe à $F_{4(-20)}$.

Ainsi on en déduit que, sous l'action de adjointe de $G$, il existe une orbite de sous-algèbres nilpotentes de dimension $1$ de $\gg$ si $g$ est hyperbolique réel, et deux orbites sinon. Dans tous les cas, l'unique orbite fermée a pour représentant une sous-algèbre de Lie de dimension $1$ incluse dans le centre de $\nn$.

\bigskip

Montrons que $\ov{\Cartan(G)}^\Ch$ est $G$-isomorphe à l'éclaté de $\bord X^{(2)}$ le long de la diagonale $\Diag(\bord X^{(2)})$. Modulo l'action de $G$, restreignons-nous au cas où l'un des points de $\bord X^{(2)}$ est le point $\xi_+ \in \bord X$ stabilisé par le sous-groupe $AN$ de $G$. Soit $\xi_- \in \bord X$ le point opposé à $\xi_+$ le long de la géodésique translatée par $A$. Alors l'application
\beq \nn & \ra & \bord X \bs \{\xi_+\} \\
Y & \mapsto & \exp(Y) \cdot \xi_- \eeq
est un homéomorphisme $AN$-équivariant. Il se prolonge en un homéomorphisme $AN$-équivariant entre la compactification de $\nn$ par $\P(\nn)$ et l'éclaté de $\bord X$ en $\xi_+$. 

Par ailleurs, on peut expliciter la compactification de Chabauty de l'espace des sous-espaces de Cartan inclus dans la sous-algèbre de Borel $\aa \oplus \nn$ grâce à l'homéomorphisme $AN$-équivariant
\beq \nn \cup \P(\nn) & \ra & \{\ll \in \ov{\Cartan(\gg)}^\Ch \,:\, \ll \subset \aa \oplus \nn\} \\
Y \in \nn & \mapsto & \Ad(\exp(Y)) \aa \\
\R (Y_\alpha + Y_{2\alpha}) \in \P(\nn) & \mapsto & \R (Y_\alpha + 2Y_{2\alpha}) \subset \nn .\eeq

Ainsi $\ov{\Cartan(G)}^\Ch$ est $G$-isomorphe à l'éclaté de $\bord X^{(2)}$ le long de $\Diag(\bord X^{(2)})$. \ep

\section{Le cas de $\SL_3(\R)$}

\label{sec:sl3}

Posons $G=\SL_3(\R)$, $K=\SO_3(\R)$ et $A$ le sous-groupe de $G$ des matrices diagonales à coefficients diagonaux strictement positifs. Soit $M$ le centralisateur de $A$ dans $K$, c'est-à-dire le sous-groupe fini de $G$ constitué des matrices diagonales à coefficients diagonaux égaux à $\pm 1$. Et soit $M'$ le normalisateur de $A$ dans $K$, c'est-à-dire le sous-groupe fini de $G$ constitué des matrices de permutation de coefficients égaux à $\pm 1$. Notons $T=MA$ le centralisateur de $A$ dans $G$, c'est un tore déployé sur $\R$ maximal, et notons $T'=M'A$ le normalisateur de $A$ (et de $T$) dans $G$.

Notons $N$ le sous-groupe de $G$ unipotent supérieur, $B=MAN$ le sous-groupe de Borel standard, c'est-à-dire le sous-groupe triangulaire supérieur, et $B_0=AN$ sa composante neutre. Notons de plus $\gg$, $\kk$, $\aa$, $\nn$ et $\bb$ les algèbres de Lie de $G$, $K$, $A$, $N$ et $B$ respectivement.

Déterminons tout d'abord les sous-algèbres abéliennes de dimension $2$ de $\bb = \aa \oplus \nn$.

\subsection{Les sous-algèbres de Lie abéliennes de $\bb$ de dimension $2$}

Soit $Y$ le sous-espace de $\Ch(\gg)$ constitué des sous-algèbres de Lie abéliennes de dimension $2$ de $\bb$. Le groupe $B_0$ agit sur $Y$ par l'action adjointe. Nous allons montrer que le sous-espace des sous-espaces de Cartan de $\gg$ inclus dans $\bb$ est dense dans $Y$, puis nous étudierons la manière dont les orbites de $B_0$ forment une structure de $CW$-complexe sur $Y$.

Notons $p_\aa : \bb = \aa \oplus \nn \ra \aa$ la projection sur $\aa$ parallèlement à $\nn$.

Notons $\alpha : \aa \ra \R$ définie par $H \mapsto H_{1,1}-H_{2,2}$ et $\beta : \aa \ra \R$ définie par $H \mapsto H_{2,2}-H_{3,3}$ : ces racines forment une base du système de racines associé au sous-espace de Cartan $\aa$. Les racines positives correspondantes sont $\Sigma^+=\{\alpha,\beta,\alpha+\beta\}$. Soit $(H_\alpha,H_\beta)$ la base de $\aa$ duale de $(\alpha,\beta)$. Leurs matrices sont $H_\alpha=\Diag(\frac{2}{3},\frac{-1}{3},\frac{-1}{3})$ et $H_\beta=\Diag(\frac{1}{3},\frac{1}{3},\frac{-2}{3})$.

Pour toute racine $\gamma \in \Sigma$, notons de plus $\aa_\gamma = \Ker \gamma$ et $A_\gamma = \exp \aa_\gamma$. Notons de plus $B'$ le sous-groupe
$$B' = \left\{ \left( \begin{array}{ccc} a & x & z \\ 0 & b & y \\ 0 & 0 & c \end{array} \right), \forall a,b,c \in \R \bs \{0\}, \forall x,y,z \in \R \,:\, \frac{a}{b} = \frac{b}{c} \right\}.$$

Notons les vecteurs
\beq
U_\alpha = \left( \begin{array}{ccc} 0 & 1 & 0 \\ 0 & 0 & 0 \\ 0 & 0 & 0 \end{array} \right), 
U_\beta = \left( \begin{array}{ccc} 0 & 0 & 0 \\ 0 & 0 & 1 \\ 0 & 0 & 0 \end{array} \right) \mbox{ et }
U_{\alpha+\beta} = \left( \begin{array}{ccc} 0 & 0 & 1 \\ 0 & 0 & 0 \\ 0 & 0 & 0 \end{array} \right).
\eeq

Notons les espaces de racines $\nn^\alpha = \R U_\alpha$, $\nn^\beta = \R U_\beta$ et $\nn^{\alpha+\beta} = \R U_{\alpha+\beta}$. Et, pour toute racine $\gamma \in \Sigma$, notons $N^\gamma = \exp \nn^\gamma$. Enfin, notons les deux sous-groupes paraboliques (maximaux à conjugaison près) :
$$ P^\alpha = \left\{\left(\begin{array}{ccc} * & * & * \\ * & * & * \\ 0 & 0 & * \end{array}\right)\right\} \mbox{ et } P^\beta = \left\{\left(\begin{array}{ccc} * & * & * \\ 0 & * & * \\ 0 & * & * \end{array}\right)\right\}.$$

Pour tout vecteur $X \in \bb$, nous noterons $X = X_\aa +x_\alpha U_\alpha + x_\beta U_\beta + x_{\alpha + \beta} U_{\alpha + \beta}$ la décomposition en espaces de racines de $X$ dans la décomposition $\bb = \aa + \nn^\alpha + \nn^\beta + \nn^{\alpha + \beta}$.

De plus, pour tout $[x:y] \in \P^1(\R)$, notons
\beq
\ll_{[x:y]} = \left\{ \left( \begin{array}{ccc} 0 & tx & s \\ 0 & 0 & ty \\ 0 & 0 & 0 \end{array} \right) \,:\, t,s \in \R \right\}. 
\eeq
On peut remarquer que les algèbres de Lie des radicaux unipotents de $P^\alpha$ et $P^\beta$ sont $\ll_{[0:1]}$ et $\ll_{[1:0]}$, ce sont des sous-algèbres de Lie de $\bb$ abéliennes de dimension $2$.

Pour toute racine $\gamma \in \Sigma^+$, notons de plus $\ll_\gamma = \aa_\gamma \oplus \nn^\gamma$.
On vérifie facilement que les algèbres de Lie $\ll_{[x:y]}$, pour $[x:y] \in \P^1(\R)$, et $\ll_\gamma$, pour $\gamma \in \Sigma^+$, appartiennent à $Y$.

\bpro \label{pro:ssalg}
Soit $\ll \in Y$. Alors il y a trois possibilités (mutuellement exclusives) :
\begin{enumerate}
\item soit il existe un unique $b \in N$ tel que $\ll = \Ad b(\aa)$ ;
\item soit il existe un unique $[x:y] \in \P^1(\R)$ tel que $\ll = \ll_{[x:y]}$ ;
\item soit il existe un unique $\gamma \in \Sigma^+$ et un unique $b \in N/N^\gamma$ tels que $\ll = \Ad b (\ll_\gamma)$.
\end{enumerate}
\epro

\bp Distinguons selon la dimension de $p_\aa(\ll)$.
\begin{enumerate}
\item Si $p_\aa(\ll)=\aa$, considérons une base de $\ll$ constituée de deux éléments diagonalisables. Alors $\ll$ est un sous-espace de Cartan de $\gg$, il existe donc $b \in G$ tel que $\ll = \Ad b (\aa)$. Puisque $\bb$ est l'algèbre de Lie d'un sous-groupe de Borel contenant $\ll$, on peut supposer de plus que $\Ad b(\bb) = \bb$, c'est-à-dire que $b \in N_G(\bb)=B=TN$. Puisque $T'$ est le normalisateur de $A$, on peut supposer que $b \in N$. Et cet élément est unique car $N \cap N_G(A) = \{e\}$.
\item Si $p_\aa(\ll)=\{0\}$, alors $\ll \subset \nn$. Puisque $\ll$ est abélienne, la projection de $\ll$ sur $\gg_\alpha + \gg_\beta$ n'est pas de dimension $2$, donc est de dimension $1$. Puisque $\ll$ ne se surjecte pas sur $\nn^\alpha + \nn^\beta$, cela implique que le noyau de la projection $\nn^{\alpha+\beta}$, qui est de dimension $1$, est inclus dans $\ll$. Il existe donc un unique $[x:y] \in \P^1(\R)$ tel que $xU_\alpha+yU_\beta \in \ll$. Ainsi $\ll=\{xtU_\alpha+ytU_\beta+zU_{\alpha+\beta} \,:\, t,z \in \R\}=\ll_{[x:y]}$.
\item Sinon $p_\aa(\ll)$ est de dimension $1$ : soit alors $X \in \ll$ tel que $p_\aa(X) = X_\aa \neq 0$. Soit $Y \in \nn$ tel que $(X,Y)$ soit une base de $\ll$. Alors, puisque $\ll$ est abélienne, la coordonnée de $[X,Y]=0$ selon $\nn^\alpha$ est $ \alpha(X_\aa) Y_\alpha = 0$ et celle selon $\nn^\beta$ est $\beta(X_\aa) Y_\beta = 0$. On se trouve alors dans l'un des trois cas suivants :
\begin{enumerate}
\item Soit $X_\aa \in \aa_\alpha$. Dans ce cas on peut supposer que $\beta(X_\aa)=1$ (c'est-à-dire $X_\aa=H_\beta$), et on a donc $Y_\beta=0$. Par ailleurs la coordonnée de $[X,Y]=0$ selon $\nn^{\alpha+\beta}$ est $y_{\alpha+\beta}-x_\beta y_\alpha=0$. Puisque $Y \neq 0$, on doit donc avoir $y_\alpha \neq 0$ : on peut supposer que $y_\alpha=1$. Quitte à ajouter à $X$ un multiple de $Y$, supposons que $x_\alpha=0$. Alors l'élément
$$ b^{-1} = \left( \begin{array}{ccc} 1 & 0 & x_{\alpha+\beta} \\ 0 & 1 & x_\beta \\ 0 & 0 & 1 \end{array} \right) \in N $$
est tel que $\Ad b^{-1}(X) = H_\beta \in \Ker \alpha$ et $\Ad b^{-1}(Y) = U_\alpha \in \nn^\alpha$. Donc $\ll = \Ad b(\aa_\alpha \oplus \nn^\alpha) = \Ad b(\ll_\alpha)$.
\item Soit $X_\aa \in \aa_\beta$. Dans ce cas symétrique au précédent, on trouve également un élément $b \in B_0$ tel que $\ll = \Ad b (\ll_\beta)$.
\item Sinon, on doit avoir $y_\alpha=y_\beta=0$ : on peut alors supposer $y_{\alpha+\beta}=1$. Puisque $X$ et $Y$ commutent, on en déduit que la coordonnée de $[X,Y]=0$ selon $\nn^{\alpha+\beta}$ est $(\alpha+\beta)(X_\aa)y_{\alpha+\beta}=(\alpha+\beta)(X_\aa)=0$.
On peut supposer que $\alpha(X_\aa)=1$ et $\beta(X_\aa)=-1$. Alors l'élément
$$ b^{-1} = \left( \begin{array}{ccc} 1 & x_\alpha & 0 \\ 0 & 1 & -x_\beta \\ 0 & 0 & 1 \end{array} \right) \in N $$
est tel que $\Ad b^{-1}(X) \in X_\aa + \nn^{\alpha+\beta}$ et $\Ad b^{-1}(Y) = Y \in \nn^{\alpha+\beta}$. Donc $\ll = \Ad b(\ll_{\alpha+\beta})$.
\end{enumerate}
Et dans chacun de ces trois cas, le normalisateur de $\ll_\gamma$ dans $N$ est égal à $N^\gamma$, d'où l'unicité de $b$ modulo $N^\gamma$. \hfill \qed
\end{enumerate}
\epp

\bcor \label{cor:representants} Sous l'action de $B_0$, il y a $8$ orbites dans $Y$, dont des représentants sont $\aa$, $\ll_{[1:0]}$, $\ll_{[0:1]}$, $\ll_{[1:1]}$, $\ll_{[1:-1]}$, $\ll_\alpha$, $\ll_\beta$ et $\ll_{\alpha+\beta}$.
\ecor

\bp Il suffit de vérifier que pour tout $[x:y] \in \P^1(\R) \bs \{[1:0],[0:1]\}$, il existe $a \in A$ tel que $\Ad a (\ll_{[x:y]}) = \ll_{[1:\pm 1]}$ (selon que $x$ et $y$ sont de même signe ou non), et que deux de ces $8$ représentants ne sont pas conjugués sous l'action adjointe de $B_0$. \ep

Considérons pour tous $x$, $y$ et $z$ réels l'élément
$$ b(x,y,z) = \exp \left( \begin{array}{ccc} 0 & x & z \\ 0 & 0 & y \\ 0 & 0 & 0 \end{array} \right) = \left( \begin{array}{ccc} 1 & x & z+\frac{xy}{2} \\ 0 & 1 & y \\ 0 & 0 & 1 \end{array} \right) \in N.$$
Le lemme suivant permet de décrire quels éléments de l'orbite de $\aa$ par $N$ convergent vers quels éléments de $Y$. 

\blem \label{lem:limite1} Soient $(x_n,y_n,z_n)_{n \in \N}$ une suite de $\R^3$ tendant vers l'infini. Pour tout $n \in \N$, notons $b_n = b(x_n,y_n,z_n)$.
\begin{enumerate}
\item Si $x_n \ra \infty$ et $(y_n,z_n+\frac{x_ny_n}{2}) \ra (y,z)$ dans $\R^2$, alors $\Ad b_n(\aa) \ral{n \ra +\infty} \Ad b(0,y,z)(\ll_\alpha)$.
\item Si $y_n \ra \infty$ et $(x_n,z_n-\frac{x_ny_n}{2}) \ra (x,z)$ dans $\R^2$, alors $\Ad b_n(\aa) \ral{n \ra +\infty} \Ad b(x,0,z)(\ll_\beta)$.
\item Si $z_n \ra \infty$ et $(x_n,y_n) \ra (x,y)$ dans $\R^2$, alors $\liml_{n \ra +\infty} \Ad b_n(\aa) \ral{n \ra +\infty} \Ad s(x,y,0)(\ll_{\alpha+\beta})$.
\item Si on n'est pas dans l'un de ces cas à extraction près, et si
$$\left[x_n(z_n+\frac{x_ny_n}{2}):y_n(-z_n+\frac{x_ny_n}{2})\right] \ra [a:b]$$
dans $\P^1(\R)$, alors $\Ad b_n(\aa) \ral{n \ra +\infty} \ll_{[a:b]}$.
\end{enumerate}
\elem

\bp Calculons les images par $\Ad b_n$ des vecteurs $H_\alpha$ et $H_\beta$ de la base de $\aa$,
\beq \Ad b_n(H_\alpha) = \left( \begin{array}{ccc} \frac{2}{3} & -x_n & -z_n+\frac{x_ny_n}{2} \\ 0 & \frac{-1}{3} & 0 \\ 0 & 0 & \frac{-1}{3} \end{array} \right) & \mbox{ et } & \Ad b_n(H_\beta) = \left( \begin{array}{ccc} \frac{1}{3} & 0 & -z_n-\frac{x_ny_n}{2} \\ 0 & \frac{1}{3} & -y_n \\ 0 & 0 & \frac{-2}{3} \end{array} \right). \eeq
\ben
\item Supposons $x_n \ra \infty$ et $(y_n,z_n+\frac{x_ny_n}{2}) \ra (y,z)$ dans $\R^2$. Les suites de vecteurs $(\Ad b_n(\frac{-1}{x_n}H_\alpha))_{n \in \N}$ et $(\Ad b_n(H_\beta))_{n \in \N}$ convergent respectivement vers $\Ad b(0,y,z)(U_\alpha)$ et $\Ad b(0,y,z)(H_\beta)$, ainsi $\liml_{n \ra +\infty} \Ad b_n(\aa) = \Ad b(0,y,z)(\ll_\alpha)$.
\item Supposons $y_n \ra \infty$ et $(x_n,z_n-\frac{x_ny_n}{2}) \ra (x,z)$ dans $\R^2$. Les suites de vecteurs $(\Ad b_n(H_\alpha))_{n \in \N}$ et $(\Ad b_n(\frac{-1}{y_n}H_\beta))_{n \in \N}$ convergent respectivement vers $\Ad b(x,0,z)(H_\alpha)$ et $\Ad b(x,0,z)(U_\beta)$, ainsi $\liml_{n \ra +\infty} \Ad b_n(\aa) = \Ad b(x,0,z)(\ll_\beta)$.
\item Supposons $z_n \ra \infty$ et $(x_n,y_n) \ra (x,y)$ dans $\R^2$. Les suites de vecteurs $(\Ad b_n(\frac{-1}{z_n}H_\alpha))_{n \in \N}$ et $(\Ad b_n(H_\alpha-H_\beta))_{n \in \N}$ convergent respectivement vers $\Ad b(x,y,0)(U_{\alpha+\beta})$ et $\Ad b(x,y,0)(H_\alpha-H_\beta)$, ainsi $\liml_{n \ra +\infty} \Ad b_n(\aa) = \Ad b(x,y,0)(\ll_{\alpha+\beta})$.
\item Sinon, considérons le vecteur $Z_n$ de $\Ad b(x_n,y_n,z_n)(\aa)$ suivant dont le coefficient $(Z_n)_{1,3}$ est nul :
\beq Z_n &=& \Ad b(x_n,y_n,z_n)\left(\left(z_n+\frac{x_ny_n}{2}\right)H_\alpha+\left(-z_n+\frac{x_ny_n}{2}\right)H_\beta \right) \\
 &=& \left( \begin{array}{ccc} \frac{z_n}{3} + \frac{x_ny_n}{2} & -x_n(z_n+\frac{x_ny_n}{2}) & 0 \\ 0 & \frac{-2z_n}{3} & -y_n(-z_n+\frac{x_ny_n}{2}) \\ 0 & 0 & \frac{z_n}{3} - \frac{x_ny_n}{2} \end{array} \right).\eeq
Montrons que, si l'on n'est pas dans l'un des trois premiers cas à extraction près, les coefficients diagonaux de $Z_n$ sont négligeables devant l'un des coefficients de la surdiagonale.
\bit
\item Si $x_n \ra \infty$, $y_n \ra y$ et $z_n+\frac{x_ny_n}{2} \ra \infty$, alors les coefficients diagonaux de $Z_n$ sont négligeables devant $(Z_n)_{1,2} = -x_n(z_n+\frac{x_ny_n}{2})$.
\item Si $x_n \ra x$, $y_n \ra \infty$ et $z_n-\frac{x_ny_n}{2} \ra \infty$, alors les coefficients diagonaux de $Z_n$ sont négligeables devant $(Z_n)_{2,3} = -y_n(-z_n+\frac{x_ny_n}{2})$.
\item Si $x_n \ra \infty$ et $y_n \ra \infty$, alors les coefficients diagonaux de $Z_n$ sont négligeables devant $(Z_n)_{1,2} = -x_n(z_n+\frac{x_ny_n}{2})$.
\eit
On peut ainsi supposer que la suite $\left(\left[x_n(z_n+\frac{x_ny_n}{2}):y_n(-z_n+\frac{x_ny_n}{2})\right]\right)_{n \in \N}$ converge vers $[a:b]$ dans $\P^1(\R)$. Dans ce cas, la suite de droites $(\R Z_n)_{n \in \N}$ converge vers $\R(a U_\alpha + b U_\beta)$ dans $\P^1(\bb)$. Par ailleurs, selon les cas, l'une des deux suites de droite $(\Ad b(x_n,y_n,z_n)(H_\alpha))_{n \in \N}$ ou $(\Ad b(x_n,y_n,z_n)(H_\alpha))_{n \in \N}$ converge vers $\R U_{\alpha+_\beta}$ dans $\P^1(\bb)$.

On en conclut que la suite de sous-algèbres de Lie $(\Ad b(x_n,y_n,z_n)(\aa))_{n \in \N}$ converge vers $\ll_{[a:b]}$. \hfill \qed
\een \epp

Le résultat suivant découle alors immédiatement du corollaire~\ref{cor:representants}.

\bcor \label{cor:orbite_dense} L'orbite de $\aa$ sous l'action adjointe de $N$ est dense dans $Y$. \hfill \qed \ecor

On en déduit alors une caractérisation intrinsèque de la compactification de Chabauty de $\Cartan(\SL_3(\R))$.

\bthm \label{thm:sgabeliens_sl3} La compactification de Chabauty $\ov{\Cartan(\SL_3(\R))}^\Ch$ coïncide avec l'espace ${\cal A}(\SL_3(\R))$. \ethm

\bp Soit $\ll$ une sous-algèbre de Lie abélienne de $\sl_3(\R)$ de dimension $2$ incluse dans une sous-algèbre de Borel. Quitte à conjuguer $\ll$ par un élément de $K$, nous pouvons supposer que $\ll$ est incluse dans la sous-algèbre de Borel standard, c'est-à-dire $\ll \in Y$. D'après le corollaire~\ref{cor:orbite_dense}, nous savons que $\ll$ est limite de sous-espaces de Cartan de $\sl_3(\R)$, donc appartient à $\ov{\Cartan(\sl_3(\R))}^\Ch$. \ep

On calcule aisément les normalisateurs dans $G$ de chacune des $8$ sous-algèbres de Lie du corollaire~\ref{cor:representants}.

\bpro \label{pro:normalisateurs}
\begin{enumerate}
\item Le normalisateur $N_G(\aa)$ de $\aa$ dans $G$ est égal à $T'=M'A$. Donc l'orbite de $\aa$ sous l'action adjointe de $B_0$ est homéomorphe à $N$, c'est-à-dire à $\R^3$.
\item Le normalisateur $N_G(\ll_{[0:1]})$ de $\ll_{[0:1]}$ dans $G$ est égal à $P^\alpha$. Donc son orbite sous l'action adjointe de $B_0$ est un point.
\item Le normalisateur $N_G(\ll_{[1:0]})$ de $\ll_{[1:0]}$ dans $G$ est égal à $P^\beta$. Donc son orbite sous l'action adjointe de $B_0$ est un point.
\item Le normalisateur $N_G(\ll_{[1:1]})$ de $\ll_{[1:1]}$ dans $G$ est égal à $B'$. Donc son orbite sous l'action adjointe de $B_0$ est homéomorphe à $\R$.
\item Le normalisateur $N_G(\ll_{[1:-1]})$ de $\ll_{[1:-1]}$ dans $G$ est égal à $B'$. Donc son orbite sous l'action adjointe de $B_0$ est homéomorphe à $\R$.
\item Pour toute racine positive $\gamma \in \Sigma^+$, le normalisateur $N_G(\ll_\gamma)$ de $\ll_\gamma$ dans $G$ est égal à $MAN^\gamma$. Donc son orbite sous l'action adjointe de $B_0$ est homéomorphe à $N/N^\gamma$, c'est-à-dire à $\R^2$.
\end{enumerate}
\epro


Le lemme suivant décrit les valeurs d'adhérence des autres orbites que $\aa$ dans $Y$ sous l'action adjointe de $B_0$, pour les représentants décrits dans le corollaire~\ref{cor:representants}.

\blem \label{lem:limite2}
\begin{enumerate}
\item Les sous-algèbres $\ll_{[0:1]}$ et $\ll_{[1:0]}$ sont normalisées par $B_0$.
\item Pour $\ll_{[1:\pm 1]}$, un système complet de représentants de $B_0$ modulo $N_{B_0}(\ll_{[1:\pm 1]}) = B' \cap B_0$ est $A_{\alpha+\beta} = \{\Diag(\lambda,\lambda^{-2},\lambda) \in A \,:\, \lambda \in \, \in \,]0,+\infty[\}$. Soit $(\lambda_n)_{n \in \N}$ une suite dans $]0,+\infty[$ et considérons pour tout $n \in \N$ l'élément $b_n = \Diag(\lambda_n,\lambda_n^{-2},\lambda_n)$ de $A_{\alpha+\beta}$. Si $\lambda_n \ra +\infty$ alors $\liml_{n \ra +\infty} \Ad h_n(\ll_{[1:\pm 1]}) = \ll_{[1:0]}$ et si $\lambda_n \ra 0$ alors $\liml_{n \ra +\infty} \Ad h_n(\ll_{[1:\pm 1]}) = \ll_{[0:1]}$.
\item Pour $\ll_\alpha$, un système complet de représentants de $B_0$ modulo $N_{B_0}(\ll_\alpha) = AN^\alpha$ est $N^\beta N^{\alpha+\beta} = \{b(0,y,z) \,:\, y,z \in \R\}$. Soit $((y_n,z_n))_{n \in \N}$ une suite de $\R^2$ tendant vers l'infini et considérons pour tout $n \in \N$ l'élément $b_n = b(0,y_n,z_n)$ de $N^\beta N^{\alpha+\beta}$. Si $[z_n:y_n^2] \ra [a:b]$ dans $\P^1(\R)$, alors $\liml_{n \ra +\infty} \Ad b(0,y_n,z_n)(\ll_\alpha) = \ll_{[a:b]}$.
\item Pour $\ll_\beta$, un système complet de représentants de $B_0$ modulo $N_{B_0}(\ll_\beta) = AN^\beta$ est $N^\alpha N^{\alpha+\beta} = \{b(x,0,z) \,:\, x,z \in \R\}$. Soit $((x_n,z_n))_{n \in \N}$ une suite de $\R^2$ tendant vers l'infini et considérons pour tout $n \in \N$ l'élément $b_n = b(x_n,0,z_n)$ de $N^\alpha N^{\alpha+\beta}$. Si $[-x_n^2:z_n] \ra [a:b]$ dans $\P^1(\R)$, alors $\liml_{n \ra +\infty} \Ad b(x_n,0,z_n)(\ll_\beta) = \ll_{[a:b]}$.
\item Pour $\ll_{\alpha+\beta}$, un système complet de représentants de $B_0$ modulo $N_{B_0}(\ll_{\alpha+\beta}) = AN^{\alpha+\beta}$ est l'ensemble $\{b(x,y,0) \,:\, x,y \in \R\}$ (ce n'est pas un sous-groupe de $B_0$). Soit $((x_n,y_n))_{n \in \N}$ une suite de $\R^2$ tendant vers l'infini et considérons pour tout $n \in \N$ l'élément $b_n = b(x_n,y_n,0)$ de $N$. Si $[-x_n:y_n] \ra [a:b]$ dans $\P^1(\R)$, alors $\liml_{n \ra +\infty} \Ad b(x_n,y_n,0)(\ll_{\alpha+\beta}) = \ll_{[a:b]}$.
\end{enumerate}
\elem

\bp
\begin{enumerate}
\item Cela a été vu dans la proposition~\ref{pro:normalisateurs}.
\item On remarque que $\Ad b_n(\ll_{[1:\pm 1]}) = \ll_{[\lambda_n^3,\pm \lambda_n^{-3}]}$, le résultat est alors clair.
\item On montre que $\R \Ad b_n\left(-y_n H_\beta + z_n U_\alpha \right) \ral{n \ra +\infty} \R (aU_\alpha + bU_\beta)$. Par ailleurs, l'une des deux suites $(\R\Ad b_n(H_\beta))_{n \in \N}$ ou $(\R\Ad b_n(U_\alpha))_{n \in \N}$ de droites de $\bb$ converge vers $\R U_{\alpha+\beta}$.
Ainsi $\liml_{n \ra +\infty} \Ad b_n(\ll_\alpha) = \ll_{[a:b]}$.
\item L'argument est le même que pour $\ll_\alpha$.
\item Le vecteur $U_{\alpha + \beta}$ est normalisé par $b_n$, et la suite $(\R \Ad b_n\left(H_\alpha -H_\beta - x_ny_n U_{\alpha + \beta} \right))_{n \in \N}$ converge vers $\R (aU_\alpha + bU_\beta)$. Ainsi $\liml_{n \ra +\infty} \Ad b_n(\ll_{\alpha+\beta}) = \ll_{[a:b]}$. \hfill \qed
\end{enumerate} \epp

\subsection{La topologie de l'espace $Y$ des sous-algèbres de Lie abéliennes de $\bb$ de dimension $2$}

Considérons l'homéomorphisme de $\R$ dans $\R$ défini par $z \mapsto z'=\sgn(z) \sqrt{|z|}$, d'inverse $z' \mapsto z=\sgn(z')(z')^2$.

Pour $n=2$ et $n=3$, notons $\B^n$ la compactification de $\R^n$ obtenue en ajoutant la sphère $\SS^{n-1}$ à l'infini (qui est la compactification géodésique de l'espace euclidien standard $\R^n$). Utilisons des coordonnées homogènes pour les sphères, en identifiant $\SS^{n-1}$ avec $(\R^n \bs \{0\})/\R_+^*$ de manière évidente. Notons $\B^2_\alpha$, $\B^2_\beta$ et $\B^2_{\alpha+\beta}$ trois copies disjointes du disque fermé $\B^2$, de disques ouverts $\R^2_\alpha$, $\R^2_\beta$ et $\R^2_{\alpha+\beta}$ et de cercles de bord $\SS^1_\alpha$, $\SS^1_\beta$ et $\SS^1_{\alpha+\beta}$. Considérons les trois applications continues
\beq f_{\alpha} : \B^3 \bs \{[\pm 1:0:0]\} & \ra & \B^2_{\alpha} \\
(x,y,z') \in \R^3 & \mapsto & \left(y,\left(z+\frac{xy}{2}\right)'\right) \in \R^2_{\alpha} \\
\ [x:y:z'] \in \SS^2 & \mapsto & \left[y:\left(z+\frac{xy}{2}\right)'\right] \in \SS^1_{\alpha},\eeq
\beq f_{\beta} : \B^3 \bs \{[0:\pm 1:0]\}  & \ra & \B^2_{\beta} \\
(x,y,z') \in \R^3  & \mapsto & \left(x,\left(z-\frac{xy}{2}\right)'\right) \in \R^2_{\beta}  \\
\ [x:y:z'] \in \SS^2  & \mapsto & \left[x:\left(z-\frac{xy}{2}\right)'\right] \in \SS^1_{\beta} \eeq
\beq \mbox{et } f_{\alpha+\beta} : \B^3 \bs \{[0:0:\pm 1]\}  & \ra & \B^2_{\alpha+\beta} \\
(x,y,z') \in \R^3 & \mapsto & (x,y) \in \R^2_{\alpha+\beta)} \\
\ [x:y:z'] \in \SS^2 & \mapsto & [x:y] \in \SS^1_{\alpha+\beta} . \eeq

\begin{figure}[h]
\begin{center}
\includegraphics[height=4cm]{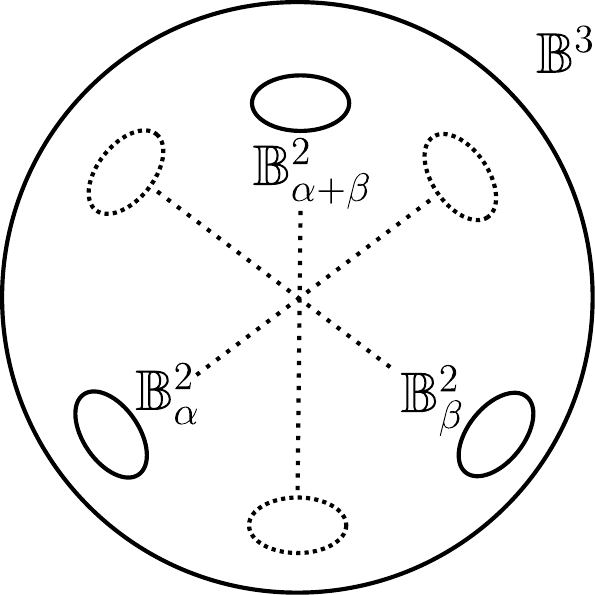}
\caption{La boule éclatée $C$.}
\label{fig:boule_eclatee}
\end{center}
\end{figure}

Remarquons que ces formules étant positivement homogènes, et les domaines de définition indiqués étant les bons, ces applications sont bien définies. Soit $C$ l'adhérence de l'image du plongement (et nous identifions l'espace de départ et l'image de ce plongement)
\beq \B^3 \bs \{[\pm 1:0:0],[0:\pm 1:0],[0:0:\pm 1]\} & \ra & \B^3 \times \B^2_\alpha \times \B^2_\beta \times \B^2_{\alpha+\beta} \\
p & \mapsto & (p,f_\alpha(p),f_\beta(p),f_{\alpha+\beta}(p)). \eeq
L'espace $C$ est la boule $\B^3$ où on a éclaté les six points $\{[\pm 1:0:0],[0:\pm 1:0],[0:0:\pm 1]\}$ en six disques, identifiés deux à deux (voir la figure~\ref{fig:boule_eclatee}).

Notons $\partial C$ le complémentaire dans $C$ de l'image de $\R^3$. Notons $C'$ le complémentaire dans $\partial C$ des trois disques $\B^2_{\alpha}$, $\B^2_{\beta}$ et $\B^2_{\alpha+\beta}$, c'est-à-dire $\SS^2$ privé des six points $[\pm 1:0:0]$, $[0:\pm 1:0]$ et $[0:0:\pm 1]$, et considérons l'application
\beq g : C' & \ra & \P^1(\R) \\
\ [x:y:z'] & \mapsto & \left[x\left(z+\frac{xy}{2}\right):y\left(-z+\frac{xy}{2}\right)\right]. \eeq
Prolongeons l'application $g$ aux bords $\SS^1_{\alpha}$, $\SS^1_{\beta}$ et $\SS^1_{\alpha+\beta}$ des plans $\R^2_{\alpha}$, $\R^2_{\beta}$ et $\R^2_{\alpha+\beta}$, par les formules suivantes.
\beq g_{\alpha} : \SS^1_{\alpha} \ra \P^1(\R) \hspace{1cm} & g_{\beta} : \SS^1_{\beta} \ra  \P^1(\R) & \hspace{1cm} g_{\alpha+\beta} : \SS^1_{\alpha+\beta} \ra \P^1(\R) \\
\ [y:z'] \mapsto [z:y^2] \hspace{1cm} & [x:z'] \mapsto [-x^2:z] & \hspace{1cm} [x:y] \mapsto [-x:y]. \eeq
%
%
On obtient ainsi une application $g$ de $C_0=C' \cup \SS^1_{\alpha} \cup \SS^1_{\beta} \cup \SS^1_{\alpha+\beta}$ dans $\P^1$.

\blem L'application $g : C_0 \ra \P^1$ est continue. \elem

\bp Il est clair que l'application $g$ est continue sur l'ouvert $C'$ de $C_0$, ainsi qu'en restriction à chaque $\SS^1_\gamma$, pour $\gamma \in \Sigma^+$.

\bigskip

Soit $[y:z']$ un point de $\SS^1_{\alpha}$, et soit $([x_n:y_n:z'_n])_{n \in \N}$ une suite de $C'$ convergeant vers $[y:z']$ dans $C$. Par continuité de $f_\beta$ et $f_{\alpha+\beta}$ en $[\pm 1:0:0]$, cela signifie que la suite $([x_n:y_n:z'_n])_{n \in \N}$ converge vers $[\pm 1:0:0]$ dans $\SS^2$, et par continuité de $f_\beta$ et $f_{\alpha+\beta}$ en $[\pm 1:0:0]$ que la suite $(f_\alpha([x_n:y_n:z'_n]))_{n \in \N}=\left(\left[y_n:\left(z_n+\frac{x_ny_n}{2}\right)'\right]\right)_{n \in \N}$ converge vers $[y:z']$ dans $\SS^1_\alpha$. Ainsi on peut supposer que la suite $(x_n)_{n \in \N}$ tend vers l'infini, que la suite $(y_n)_{n \in \N}$ converge vers $y$, et que la suite $(z_n+\frac{x_ny_n}{2})_{n \in \N}$ converge vers $z$. Or $g([x_n:y_n:z'_n])$ est égal à $\left[z_n+\frac{x_ny_n}{2}:y_n \left(-\frac{z_n}{x_n}+\frac{y_n}{2}\right)\right]$.

Si $y \neq 0$, alors la suite $(z_n)_{n \in \N}$ est équivalente à $(-\frac{x_ny}{2})_{n \in \N}$, donc la suite $(-\frac{z_n}{x_n}+\frac{y_n}{2})_{n \in \N}$ converge vers $y$, et donc la suite $(g([x_n:y_n:z'_n]))_{n \in \N}$ converge vers $[z:y^2]=g_\alpha([y:z'])$.
Si $y=0$, alors on en déduit que la suite $(z_n)_{n \in \N}$ est négligeable devant $(x_n)_{n \in \N}$, et donc la suite $(g([x_n:y_n:z'_n]))_{n \in \N}$ converge vers $[1:0]=g_{\alpha}([y:z'])$.

\bigskip

Pour $\SS^1_{\beta}$, la continuité de $g$ se montre de manière identique.

\bigskip

Soit $[x:y]$ un point de $\SS^1_{\alpha+\beta}$, et soit $([x_n:y_n:z'_n])_{n \in \N}$ une suite de $C'$ convergeant vers $[x:y]$ dans $C$. Cela signifie que la suite $([x_n:y_n:z'_n])_{n \in \N}$ converge vers $[0:0:\pm 1]$ dans $\SS^2$, et que la suite $(f_{\alpha+\beta}([x_n:y_n:z'_n]))_{n \in \N} = ([x_n:y_n])_{n \in \N}$ converge vers $[x:y]$ dans $\SS^1$. On peut donc supposer que la suite $(x_n)_{n \in \N}$ converge vers $x$, que la suite $(y_n)_{n \in \N}$ converge vers $y$, et que la suite $(z_n)_{n \in \N}$ tend vers l'infini. Or $g([x_n:y_n:z'_n])$ est égal à $\left[x_n \left(z_n+\frac{x_ny_n}{2}\right):y_n \left(-z_n+\frac{x_ny_n}{2}\right)\right]$. Puisque la suite $(x_ny_n)_{n \in \N}$ est négligeable devant la suite $(z_n)_{n \in \N}$, on en déduit que la suite $(g([x_n:y_n:z'_n]))_{n \in \N}$ est équivalente à la suite $([x_n:-y_n])_{n \in \N}$, donc converge vers $[x:-y]=g_{\alpha+\beta}([x:y])$.

\bigskip

On a donc montré que l'application $g$ était continue sur $C_0$.
\ep

Considérons l'application suivante.
\beq \phi : C & \ra & Y \\
(x,y,z') \in \R^3 & \mapsto & \Ad b(x,y,z) \aa \\
(y,z') \in \R^2_{\alpha} & \mapsto & \Ad b(0,y,z) \ll_\alpha \\
(x,z') \in \R^2_{\beta} & \mapsto & \Ad b(x,0,z) \ll_\beta \\
(x,y) \in \R^2_{\alpha+\beta} & \mapsto & \Ad b(x,y,0) \ll_{\alpha+\beta} \\
c \in C_0 & \mapsto & \ll_{g(c)}. \eeq

\bthm L'application $\phi$ est continue, surjective et induit un homéomorphisme $\widetilde{\phi}$ du recollement $C \cup_g \P^1(\R)$ sur l'espace topologique $Y$. \ethm

\bp Puisque l'application $\P^1(\R) \ra Y$ définie par $[a:b] \mapsto \ll_{[a:b]}$ est un plongement, et que l'application $\phi$ factorise par $g$ sur $C_0$, on en déduit que $\phi$ est continue en restriction à $C_0$. L'application $\phi$ est également continue sur l'ouvert $\R^3$ de $C$.
\beq \B^2_{\alpha} & \ra & Y \\
(y,z') \in \R^2_{\alpha} & \mapsto & \Ad b(0,y,z')\ll_\alpha\\
\ [y:z'] \in \SS^1_{\alpha} & \mapsto & \ll_{g_\alpha([y:z'])}, \\
\B^2_{\beta} & \ra & Y \\
(x,z') \in \R^2_{\beta} & \mapsto & \Ad b(x,0,z')\ll_\beta\\
\ [x:z'] \in \SS^1_{\beta} & \mapsto & \ll_{g_\beta([x:z'])} \\
\mbox{et } \hspace{1cm} \B^2_{\alpha+\beta} & \ra & Y \\
(x,y) \in \R^2_{\alpha+\beta} & \mapsto & \Ad b(x,y,0)\ll_{\alpha+\beta}\\
\ [x:y] \in \SS^1_{\alpha+\beta} & \mapsto & \ll_{g_{\alpha+\beta}([x:y])}.  \eeq
D'après le lemme~\ref{lem:limite2}, ce sont des plongements. Puisqu'en restriction aux trois disques $\B^2_{\alpha}$, $\B^2_{\beta}$ et $\B^2_{\alpha+\beta}$, l'application $\phi$ factorise par ces plongements, on en déduit que l'application $\phi$ est continue en restriction à $C \bs \R^3$.

\bigskip

Soit $(y,z')$ un point de $\R^2_{\alpha}$, et soit $(x_n,y_n,z'_n)_{n \in \N}$ une suite de $\R^3$ convergeant vers $(y,z')$ dans $C$. Cela signifie que la suite $(x_n,y_n,z'_n)_{n \in \N}$ converge vers $[\pm 1:0:0]$ dans $\B^3$, et que la suite $(f_\alpha(x_n,y_n,z'_n))_{n \in \N}$ converge vers $(y,z')$ dans $\R^2_{\alpha}$. Puisque $(x_n)_{n \in \N}$ tend vers l'infini et que la suite $\left(y_n,z_n+\frac{x_ny_n}{2}\right)_{n \in \N}$ converge vers $(y,z)$, d'après le lemme~\ref{lem:limite1}(1), on en déduit que la suite de sous-algèbres de Lie$(\phi(x_n,y_n,z_n))_{n \in \N}$ converge vers $\Ad s(0,y,z) \ll_\alpha = \phi(y,z')$.

\bigskip

Soit $[y:z']$ un point de $\SS^1_{\alpha}$, et soit $(x_n,y_n,z'_n)_{n \in \N}$ une suite de $\R^3$ convergeant vers $[y:z']$ dans $C$. Cela signifie que la suite $(x_n,y_n,z'_n)_{n \in \N}$ converge vers $[\pm 1:0:0]$ dans $\B^3$, et que la suite $(f_\alpha(x_n,y_n,z'_n))_{n \in \N}$ converge vers $[y:z']$ dans $\B^2_{\alpha}$. 

Supposons $y \neq 0$, alors la suite $\left(\frac{z_n+\frac{x_ny_n}{2}}{y_n^2}\right)_{n \in \N}$ converge vers $\frac{z}{y^2}$, or la suite $(\frac{x_n}{y_n})_{n \in \N}$ tend vers l'infini, donc la suite $(\frac{z_n}{y_n^2})_{n \in \N}$ est équivalente à $(-\frac{x_n}{2y_n})_{n \in \N}$. Ainsi dans $\P^1(\R)$ la suite $\left(\left[x_n(z_n+\frac{x_ny_n}{2}):y_n(-z_n+\frac{x_ny_n}{2})\right]\right)_{n \in \N}$ est équivalente à la suite $\left(\left[x_n \frac{z}{y^2}:y_n\frac{x_n}{y_n}\right]\right)_{n \in \N}$, donc converge vers $[z:y^2] = g_{\alpha}([y:z'])$. D'après le lemme~\ref{lem:limite1}(4), on en déduit que la suite de sous-algèbres de Lie $(\phi(x_n,y_n,z_n))_{n \in \N}$ converge vers $\ll_{[z:y^2]} = \phi([y:z'])$.

Supposons $y=0$, alors la suite $\left(\frac{y_n^2}{z_n+\frac{x_ny_n}{2}}\right)_{n \in \N}$ converge vers $0$. On en déduit que la suite $\left(\left[x_n(z_n+\frac{x_ny_n}{2}):y_n(-z_n+\frac{x_ny_n}{2})\right]\right)_{n \in \N} = \left(\left[x_n(z_n+\frac{x_ny_n}{2}):-y_n(z_n+\frac{x_ny_n}{2})+x_ny_n^2\right]\right)_{n \in \N}$
converge vers $[1:0]=g_{\alpha}([y:z'])$. D'après le lemme~\ref{lem:limite1}(4), on en déduit que la suite de sous-algèbres de Lie $(\phi(x_n,y_n,z_n))_{n \in \N}$ converge vers $\ll_{[1:0]} = \phi([y:z'])$.

\bigskip

La continuité de $\phi$ sur $\B^2_{\beta}$ est identique.

\bigskip

Soit $(x,y)$ un point de $\R^2_{\alpha+\beta}$, et soit $(x_n,y_n,z'_n)_{n \in \N}$ une suite de $\R^3$ convergeant vers $(x,y)$ dans $C$. Cela signifie que la suite $(x_n,y_n,z'_n)_{n \in \N}$ converge vers $[0:0:\pm 1]$ dans $\B^3$, et que la suite $(f_{\alpha+\beta}(x_n,y_n,z'_n))_{n \in \N}=(x_n,y_n)_{n \in \N}$ converge vers $(x,y)$ dans $\R^2_{\alpha+\beta}$. Puisque $(z_n)_{n \in \N}$ tend vers l'infini et comme la suite $(x_n,y_n)_{n \in \N}$ converge vers $(x,y)$, d'après le lemme~\ref{lem:limite1}(3), on en déduit que la suite de sous-algèbres de Lie $(\phi(x_n,y_n,z_n))_{n \in \N}$ converge vers $\Ad b(x,y,0)(\ll_{\alpha+\beta}) = \phi(x,y)$.

\bigskip

Soit $[x:y]$ un point de $\SS^1_{\alpha+\beta}$, et soit $(x_n,y_n,z'_n)_{n \in \N}$ une suite de $\R^3$ convergeant vers $(x,y)$ dans $C$. Cela signifie que la suite $(x_n,y_n,z'_n)_{n \in \N}$ converge vers $[0:0:\pm 1]$ dans $\B^3$, et que la suite $(f_{\alpha+\beta}(x_n,y_n,z'_n))_{n \in \N}=(x_n,y_n)_{n \in \N}$ converge vers $[x:y]$ dans $\B^2_{\alpha+\beta}$. Alors la suite $$\left(\left[x_n(z_n+\frac{x_ny_n}{2}):y_n(-z_n+\frac{x_ny_n}{2})\right]\right)_{n \in \N}$$
est équivalente à la suite $([x_n:-y_n])_{n \in \N}$, donc converge vers $[x:-y]=g_{\alpha+\beta}([x:y])$ dans $\P^1(\R)$. Donc d'après le lemme~\ref{lem:limite1}(4), on en déduit que la suite de sous-algèbres de Lie $(\phi(x_n,y_n,z_n))_{n \in \N}$ converge vers $\ll_{[x:-y]} = \phi([x:y])$.

\bigskip

On a donc montré que l'application $\phi$ était continue sur $C$. Puisqu'elle factorise par $g$ sur $C_0$, on en déduit qu'elle induit une application continue $\widetilde{\phi}$ de $C \cup_g \P^1(\R)$ sur $Y$.

L'application $\phi$ est une bijection de $\R^3$ sur l'orbite de $\aa$ dans $Y$ sous l'action adjointe de $B_0$. Pour tout $\gamma \in \Sigma^+$, l'application $\phi$ est une bijection de $\R^2_\gamma$ sur l'orbite de $\ll_\gamma$ dans $Y$ sous l'action adjointe de $B_0$. Et l'application $\phi$ est une bijection de $\P^1(\R)$ sur le sous-espace $\{\ll_{[a:b]} \,:\, [a:b] \in \P^1(\R)\}$ de $Y$. Ainsi d'après la proposition~\ref{pro:ssalg}, l'application $\phi$ est bijective. Or l'espace $C \cup_g \P^1(\R)$ est compact et l'espace $Y$ est séparé, donc l'application $\widetilde{\phi}$ est un homéomorphisme.
\ep

\vfill

\begin{figure}[h]
\begin{center}
\includegraphics[height=5.5cm]{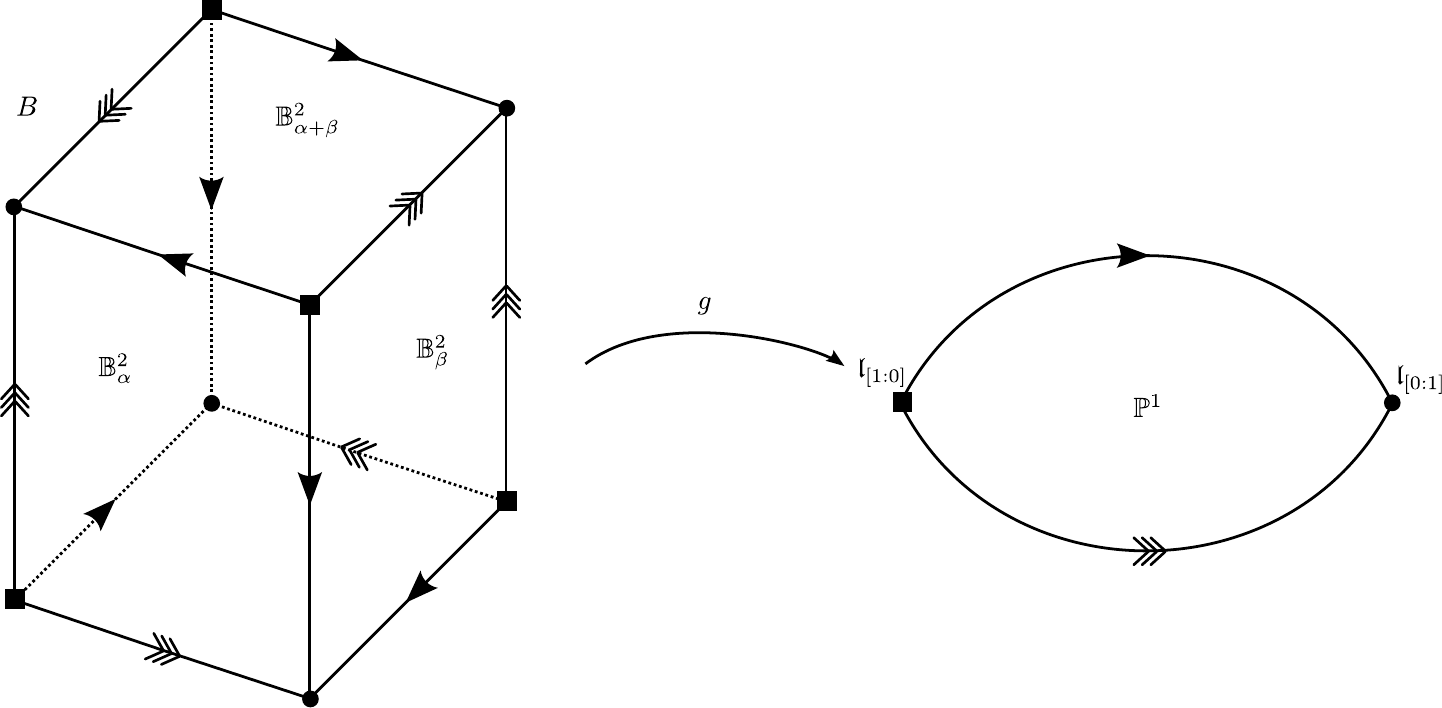}
\caption{Le $CW$-complexe $Y$}
\label{fig:cwy}
\end{center}
\end{figure}

\bcor \label{cor:CWcomplexe} L'espace $Y$ est un CW-complexe (voir la figure~\ref{fig:cwy}), dont les simplexes sont donnés par les $B_0$-orbites. Les voici, classés par dimensions.
\begin{itemize}
\item Il y a deux $0$-simplexes, qui sont les points $\ll_{[1:0]}$ et $\ll_{[0:1]}$.
\item Il y a deux $1$-simplexes, qui sont les orbites de $\ll_{[1:\pm 1]}$.
\item Il y a trois $2$-simplexes, qui sont les orbites de $\Ker \gamma \oplus \gg_\gamma$ pour $\gamma \in \Sigma^+$.
\item Il y a un $3$-simplexe, qui est l'orbite de $\aa$.
\end{itemize} \hfill \qed
\ecor

\bpro \label{pro:pi1y} Le $2$-squelette $Y^{(2)}$ de $Y$ a le type d'homotopie qu'un bouquet formé de deux sphères $\SS^2$ et d'un plan projectif réel $\P^2(\R)$. En particulier, l'espace $Y$ a pour groupe fondamental $\pi_1(Y) = \Z / 2\Z$, dont un générateur est donné par la classe d'homotopie d'un lacet parcourant une fois le cercle $\P^1(\R)$. \epro

\bp Le $2$-squelette $Y^{(2)}$ de $Y$ est le recollement des trois $2$-cellules $\R^2_\alpha$, $\R^2_\beta$ et $\R^2_{\alpha+\beta}$ par les trois applications $g_\alpha$, $g_\beta$ et $g_{\alpha+\beta}$ sur le $1$-squelette $Y^{(1)} \simeq \P^1(\R)$.

Or les applications $g_\alpha$ et $g_\beta$ sont homotopes à des applications constantes, et le recollement de $\R^2_{\alpha+\beta}$ sur $\P^1(\R)$ par $g_{\alpha+\beta}$ est homéomorphe à $\P^2(\R)$. Ainsi le $2$-squelette $Y^{(2)}$ a le type d'homotopie d'un bouquet formé de deux sphères $\SS^2$ et d'un plan projectif $\P^2(\R)$. Ainsi le groupe fondamental de $Y$ est isomorphe au groupe fondamental de $\P^2(\R)$, c'est-à-dire $\Z / 2\Z$. \ep 
%
%
%
%

\subsection{Compactification de $G/T'$}

Notons $Z = \ov{\Cartan(\SL_3(\R))}^\Ch$ la compactification de Chabauty de $G/T' = \Cartan(\SL_3(\R))$. Considérons l'application continue
\beq \pi : G \times Y & \ra & Z \\
(g,y) & \mapsto & g \cdot y, \eeq
où le groupe $G$ agit sur $Z$ par l'action adjointe. Puisque toute $G$-orbite dans $Z$ rencontre $Y$, l'application $\pi$ est surjective et ouverte. Choisissons comme point base $(e,y_0)=(e,\aa)$ dans $G \times Y$ et $z_0=\pi(e,y_0)=\aa$ dans $Z$.

%

Commençons par démontrer un lemme de topologie générale.

\blem \label{lem:pi_surjective_top} Soient $(E,e_0)$ est un espace topologique pointé localement connexe par arcs, $(F,f_0)$ un espace topologique pointé localement simplement connexe et $\pi : (E,e_0) \ra (F,f_0)$ une application continue, surjective et ouverte, telle que $\pi^{-1}(f_0)$ soit connexe par arcs. Alors l'application $\pi$ induit une surjection
$$ \pi_* : \pi_1(E,e_0) \ra \pi_1(F,f_0).$$
\elem

\bp Montrons qu'on peut relever localement les chemins, à homotopie près. Plus précisément, montrons que, pour tout $e \in E$, il existe un voisinage $U$ de $e$ dans $E$ tel que, pour tout $e' \in U$ et tout chemin $c$ de $\pi(e)$ à $\pi(e')$ dans $\pi(U)$, il existe un chemin $\widetilde{c}$ de $e$ à $e'$ dans $U$ tel que $\pi(\widetilde{c})$ soit homotope dans $F$ à $c$ relativement aux extrémités.

\bigskip

Fixons $e \in E$ et $f=\pi(e)$. L'espace $F$ est localement simplement connexe : soit donc un voisinage simplement connexe $V$ de $f$ dans $F$. L'espace $E$ est localement connexe par arcs : soit donc un voisinage connexe par arcs $U$ de $e$ dans $E$ tel que $\pi(U) \subset V$.

\bigskip

Soient $e' \in U$ et $f'=\pi(e')$. Soit $c:[0,1] \ra \pi(U)$ un chemin continu de $f$ à $f'$. Considérons un chemin continu $\widetilde{c} : [0,1] \ra U$ de $e$ à $e'$. Alors $\pi(\widetilde{c})$ et $c$ sont deux chemins continus de $f$ à $f'$ inclus dans $\pi(U)$, donc inclus dans $V$ : ils sont ainsi homotopes dans $F$ relativement à leurs extrémités.

\bigskip

Ainsi, pour tout lacet $\ell$ de $F$ basé en $f_0$, il existe un chemin $\widetilde{c}$ de $E$ d'origine $e_0$ et d'extrémité appartenant à $\pi^{-1}(f_0)$ tel que les deux lacets $\ell$ et $\pi(\widetilde{c})$ (basés en $f_0$) soient homotopes dans $F$. Or $\pi^{-1}(f_0)$ est connexe par arcs : soit $\widetilde{c'}$ un chemin dans $\pi^{-1}(f_0)$ joignant les deux extrémités de $\widetilde{c}$, alors le lacet concaténé $\widetilde{\ell} = \widetilde{c'} \cdot \widetilde{c}$ basé en $e_0$ est tel que $\pi(\widetilde{\ell})$ soit homotope à $\ell$ dans $F$. L'application $\pi_*$ est donc surjective.
\ep

\bcor \label{cor:pi_surjective} L'application $\pi : (G \times Y,(e,y_0)) \ra (Z,z_0)$ induit une surjection
$$\pi_* : \pi_1(G \times Y,(e,y_0)) \ra \pi_1(Z,z_0).$$
\ecor

\bp Il suffit de vérifier les hypothèses du lemme~\ref{lem:pi_surjective_top} pour l'application $\pi : (G \times Y,(e,y_0)) \ra (Z,z_0)$. Les espaces $Y$ et $Z$ sont des variétés semi-algébriques réelles, donc sont localement contractiles. De plus, la préimage $\pi^{-1}(z_0)=\{(g,h \cdot \aa) \in G \times Y \,:\, gh \in T'\} \simeq G \times G/T'$ est connexe par arcs. Ainsi d'après le lemme~\ref{lem:pi_surjective_top} l'application $\pi_* : \pi_1(G \times Y,(e,y_0)) \ra \pi_1(Z,z_0)$ est surjective. \ep

\bthm \label{thm:chabautysl3} La compactification de Chabauty $Z=\ov{\Cartan(\SL_3(\R))}^\Ch$ de $G/T' = \Cartan(\SL_3(\R))$ est simplement connexe. \ethm

\bp
Il suffit de montrer que les générateurs de $\pi_1(G) \simeq \Z/2\Z$ et de $\pi_1(Y) \simeq \Z/2\Z$ ont une image triviale par $\pi_*$ dans $\pi_1(Z)$. Chacun de ces générateurs est la classe d'homotopie d'un lacet dans $G \times Y^{(1)}$. Rappelons que $Y^{(1)}$ est le $1$-squelette de $Y$, c'est-à-dire $\P^1(\R)$.

\bigskip

D'après la proposition~\ref{pro:normalisateurs}, l'image par $\pi$ du sous-espace $G \times Y^{(1)}$ de $G \times Y$ est homéomorphe au quotient $Q=(G/B' \times \P^1(\R))/\!\sim$, avec les identifications suivantes.
\beq (gB',[0:1]) \sim (g'B',[0:1])  \mbox{ si }  gP^\alpha = g'P^\alpha, & (gB',[1:0])  \sim  (g'B',[1:0]) \mbox{ si }  gP^\beta = g'P^\beta \\
\mbox{et }\forall [s:t] \in \P^1(\R) \bs \{[0:1],[1:0]\}, & (gB',[s:t])  \sim  (g'B',[-s:t]) \mbox{ si } gB' = \sigma g'B', \eeq
où $\sigma \in M \bs M_{\alpha = \beta}$ est tel que $\Ad \sigma (\ll_{[1:1]}) = \ll_{[-1:1]}$, par exemple
$$ \sigma = \left( \begin{array}{ccc} -1 & 0 & 0 \\ 0 & -1 & 0 \\ 0 & 0 & 1 \end{array} \right) \in M.$$

Montrons qu'un générateur de $\pi_1(Y)$ a une image par $\pi_*$ triviale dans $\pi_1(Z)$. Considérons le lacet $c$ suivant dans $Q$, dont la classe d'homotopie engendre $\pi_*(\pi_1(Y))$, effectuant une fois le tour du cercle $Y^{(1)} \simeq \P^1(\R)$, composé des deux chemins $t \in [0,1] \mapsto \left[B',[t:1-t]\right]$ et $t \in [0,1] \mapsto \left[B',[1-t:-t]\right]$. D'après les équivalences définissant $Q$, le lacet $c$ est aussi la composée des deux chemins $t \in [0,1] \mapsto \left[B',[t:1-t]\right]$ et $t \in [0,1] \mapsto \left[\sigma B',[1-t:t]\right]$.

Fixons un chemin $c'$ de $B'$ à $\sigma B'$ dans $G/B'$. Alors le lacet $c$ est homotope à la composée des quatre chemins
\beq t \in [0,1] \mapsto \left[B',[t:1-t]\right], && t \in [0,1] \mapsto \left[c'(t),[1:0]\right], \\
t \in [0,1] \mapsto \left[\sigma B',[1-t:t]\right] & \mbox{ et } & t \in [0,1] \mapsto \left[c'(1-t),[0:1]\right].\eeq
Or la composée des deuxième et troisième chemins est homotope à la composée des chemins $t \in [0,1] \mapsto \left[B',[1-t:t]\right]$ et $t \in [0,1] \mapsto \left[c'(t),[0:1]\right]$. Ainsi le lacet $c$ est homotope à zéro.

\bigskip

Montrons qu'un générateur de $\pi_1(G)$ a une image par $\pi_*$ triviale dans $\pi_1(Z)$. Considérons le lacet $d$ suivant dans $Q$, dont la classe d'homotopie engendre $\pi_*(\pi_1(G))$ : 
\beq d : [0,1] & \ra & Q \\
t & \mapsto & \left[\left(\begin{array}{ccc} 1&0&0 \\ 0& \cos 2\pi t & -\sin 2\pi t \\ 0 & \sin 2\pi t & \cos 2\pi t \end{array}\right)B',[0:1]\right].\eeq
Ce lacet est homotope à la composée des trois chemins
\beq 
t \in [0,1] & \mapsto & \left[B',[t:1-t]\right],\\
t \in [0,1] & \mapsto & \left[\left(\begin{array}{ccc} 1&0&0 \\ 0& \cos 2\pi t & -\sin 2\pi t \\ 0 & \sin 2\pi t & \cos 2\pi t \end{array}\right)B',[1:0]\right] \\
\mbox{et } t \in [0,1] & \mapsto & \left[B',[1-t:t]\right].\eeq
Le deuxième chemin est constant dans $Q$, ainsi le lacet $d$ est homotope à zéro.

\bigskip

D'après le corollaire~\ref{cor:pi_surjective}, l'application $\pi_*$ est surjective, et nous venons de montrer que son image était triviale, donc l'espace $Z$ est simplement connexe.
\ep

Ce résultat peut sembler surprenant, étant donné que l'espace $\Plats(X)$ n'est pas simplement connexe : son groupe fondamental est de cardinal $48$ (voir le lemme~\ref{lem:pi1_plats}).

\section{Le cas de $\SL_4(\R)$}

\label{sec:sl4}


Posons $G=\SL_4(\R)$, $K=\SO_4(\R)$ et $A$ le sous-groupe de $G$ des matrices diagonales à coefficients diagonaux strictement positifs. Soit $M$ le centralisateur de $A$ dans $K$, c'est-à-dire le sous-groupe fini de $G$ constitué des matrices diagonales à coefficients diagonaux égaux à $\pm 1$.

Notons $N$ le sous-groupe de $G$ unipotent supérieur, $B=MAN$ le sous-groupe de Borel standard, c'est-à-dire le sous-groupe triangulaire supérieur, et $B_0=AN$ sa composante neutre.

Notons de plus $\gg$, $\kk$, $\aa$, $\nn$ et $\bb$ les algèbres de Lie de $G$, $K$, $A$, $N$ et $B$ respectivement.

Tout d'abord, remarquons qu'il existe une sous-algèbre de Lie de $\gg$ abélienne maximale incluse dans $\bb$ de dimension $4$ :
$$ \ll_0 = \left\{ \left( \begin{array}{cccc} 0 & 0 & * & * \\ 0 & 0 & * & * \\ 0&0&0&0 \\ 0&0&0&0 \end{array} \right) \right\}.$$
Puisque la dimension des sous-espaces de Cartan est $3$, les éléments de $\ov{\Cartan(\sl_4(\R))}^\Ch$ ne sont pas toujours des sous-algèbres de Lie abéliennes maximales, contrairement au cas de $\sl_3(\R)$. D'après~\cite{winternitz_zassenhaus}, $\ll_0$ est (à conjugaison près) la seule sous-algèbre abélienne maximale de $\sl_4(\R)$ de dimension $4$. Nous allons étudier les sous-algèbres de Lie de $\gg$ abéliennes incluses dans $\bb$ de dimension $3$.

Notons $p_\aa : \bb = \aa \oplus \nn \ra \aa$ la projection sur $\aa$ parallèlement à $\nn$.

Considérons les trois racines de $\aa$ définies par $\alpha(H) = H_{1,1}-H_{2,2}$, $\beta(H)= H_{2,2}-H_{3,3}$ et $\gamma(H) = H_{3,3}-H_{4,4}$, elles forment une base $\Delta$ du système de racines $\Sigma$ associé au sous-espace de Cartan $\aa$. Les racines positives correspondantes sont $\Sigma^+=\{\alpha,\beta,\gamma,\alpha+\beta,\beta+\gamma,\alpha+\beta+\gamma\}$. Pour toute racine positive $\delta$, notons $\nn^\delta$ l'espace de racine de $\delta$ et notons $U_\delta \in \nn^\delta$ la matrice ayant un coefficient égal à $1$ en position $\delta$, et $0$ ailleurs.

Pour toute partie $I \subset \Delta$, notons de plus $\aa_I = \cap_{\delta \in I} \Ker \delta$ et $\aa^I = \oplus_{\delta \in I} \R H_\delta$ (l'orthogonal de $\aa_I$ dans $\aa$ pour la forme de Killing), ainsi que $A_I = \exp \aa_I$ et $A^I = \exp(\aa^I)$. Notons $\Sigma^{I,+}$ l'ensemble des racines positives s'écrivant comme somme de racines de $I$. Notons $\nn^I = \oplus_{\delta \in \Sigma^{I,+}} \nn^\delta$.  


Notons, pour toute racine $\delta \in \Sigma$, l'espace de racine $\nn^\delta = \R U_\delta$ et le sous-groupe $N^\delta = \exp \nn^\delta$. Soit $I \subset \Delta$, notons $\Sigma^{I,+}$ l'ensemble des racines positives s'écrivant comme somme de racines de $I$. Notons $\nn^I = \oplus_{\delta \in \Sigma^{I,+}} \nn^\delta$ et $\nn_I = \oplus_{\delta \in \Sigma^+ \bs \Sigma^{I,+}} \nn^\delta$.

Nous considérerons toujours les projections orthogonales par rapport à la forme de Killing. En particulier, notons $p_{\alpha,\gamma} : \bb \ra \nn^\alpha \oplus \nn^\gamma$ la projection sur $\nn^\alpha \oplus \nn^\gamma$ parallèlement à $\aa \oplus \bigoplus_{\delta \in \Sigma^+ \bs \{\alpha,\gamma\}} \nn^\delta$. Et, pour toute racine positive $\delta \in \Sigma^+$, notons $p_\delta : \bb \ra \nn^\delta$ la projection sur $\nn^\delta$ parallèlement à $\aa \oplus \nn_\delta$.

Définissons, pour tout $[x:y:z] \in \P^2(\R)$ avec $x$ et $z$ non nuls, la sous-algèbre de Lie de $\gg$
$$ \ll_{[x:y:z]} = \left\{ \left( \begin{array}{cccc} 0 & ax & bx & c \\ 0 & 0 & ay & bz \\ 0&0&0&az \\ 0&0&0&0 \end{array} \right) \,:\, a,b,c \in \R \right\}.$$
Et définissons, pour tout $[x:y:z:t] \in \P^3(\R)$, la sous-algèbre de Lie de $\gg$
$$ \ll_{[x:y:z:t]} = \left\{ \left( \begin{array}{cccc} 0 & 0 & a & b \\ 0 & 0 & c & d \\ 0&0&0&0 \\ 0&0&0&0 \end{array} \right) \,:\, a,b,c,d \in \R \,:\, ax+by+cz+dt=0 \right\}.$$
Définissons également, pour tout $(y,t) \in \R^2$, les sous-algèbres de Lie de $\gg$
$$ \ll_{\alpha,y,t} = \left\{ \left( \begin{array}{cccc} 0 & a & b & c \\ 0 & 0 & ay & at \\ 0&0&0&0 \\ 0&0&0&0 \end{array} \right) \,:\, a,b,c \in \R \right\} \mbox{ et } \ll_{\gamma,y,t} = \left\{ \left( \begin{array}{cccc} 0 & 0 & at & c \\ 0 & 0 & ay & b \\ 0&0&0&a \\ 0&0&0&0 \end{array} \right) \,:\, a,b,c \in \R \right\}. $$
Et, pour tout $(x,y) \in \R^2$, définissons la sous-algèbre de Lie de $\gg$
$$ \ll_{x,y} = \left\{ \left( \begin{array}{cccc} 0 & a & by & c \\ 0 & 0 & 0 & ax \\ 0&0&0&b \\ 0&0&0&0 \end{array} \right) \,:\, a,b,c \in \R \right\}.$$
Ce sont toutes des sous-algèbres de Lie de $\gg$ abéliennes de dimension $3$ incluses dans $\bb$.

Pour toute paire de racines primitives distinctes $\{\delta,\delta'\} \subset \Delta$, notons $\zz^{\delta,\delta'}$ l'algèbre de Lie dérivée du centralisateur de $\aa_{\Delta \bs \{\delta,\delta'\}}$ dans $\gg$ :
\beq \zz^{\beta,\gamma} &=& \left\{ \left( \begin{array}{cccc} 0&0&0&0 \\ 0&*&*&* \\ 0&*&*&* \\ 0&*&*&* \end{array} \right) \right\} \simeq \sl_3(\R),\\
\zz^{\alpha,\gamma} &=& \left\{ \left( \begin{array}{cccc} *&*&0&0 \\ *&*&0&0 \\ 0&0&*&* \\ 0&0&*&* \end{array} \right) \right\} \simeq \sl_2(\R) \oplus \sl_2(\R) \\
\mbox{et } \zz^{\alpha,\beta} &=& \left\{ \left( \begin{array}{cccc} *&*&*&0 \\ *&*&*&0 \\ *&*&*&0 \\ 0&0&0&0 \end{array} \right) \right\} \simeq \sl_3(\R),\eeq
et notons $Z^{\beta,\gamma}$, $Z^{\alpha,\gamma}$ et $Z^{\alpha,\beta}$ les trois sous-groupes de Lie connexes de $G$ d'algèbres de Lie $\zz^{\beta,\gamma}$, $\zz^{\alpha,\gamma}$ et $\zz^{\alpha,\beta}$ respectivement.

Définissons de plus, pour tout $[x:y] \in \P^1(\R)$, les trois sous-algèbres de Lie de $\gg$ abéliennes de dimension $2$
\beq \ll_{[x:y]}^{\alpha,\beta} &=& \left\{ \left( \begin{array}{cccc} 0 & ax & b & 0 \\ 0 & 0 & ay & 0 \\ 0&0&0&0 \\ 0&0&0&0 \end{array} \right) \,:\, a,b \in \R \right\} \subset \zz^{\alpha,\beta} \cap \nn, \\
 \ll^{\alpha,\gamma} &=& \left\{ \left( \begin{array}{cccc} 0 & a & 0 & 0 \\ 0 & 0 & 0 & 0 \\ 0&0&0&b \\ 0&0&0&0 \end{array} \right) \,:\, a,b \in \R \right\} \subset \zz^{\alpha,\gamma} \cap \nn\\
 \mbox{et } \ll_{[x:y]}^{\beta,\gamma} &=& \left\{ \left( \begin{array}{cccc} 0 & 0 & 0 & 0 \\ 0 & 0 & ax & b \\ 0&0&0&ay \\ 0&0&0&0 \end{array} \right) \,:\, a,b \in \R \right\} \subset \zz^{\beta,\gamma} \cap \nn.\eeq

\bpro \label{pro:ssalg_sl4}
Soit $\ll$ une sous-algèbre de Lie de $\gg$ abélienne de dimension $3$ incluse dans $\bb$. Alors il y a dix possibilités (mutuellement exclusives) :
\begin{enumerate}
\item soit il existe un unique $b \in N$ tel que $\ll = \Ad b(\aa)$ ;
\item soit il existe une unique racine primitive $\delta \in \Delta$ et un unique $b \in N_\delta$ tels que $\ll = \Ad b(\aa_\delta \oplus \nn^\delta)$ ;
\item soit il existe un unique $[x:y] \in \P^1(\R)$ et un unique $b \in N_{\alpha,\beta}$ tels que $$\ll = \Ad b(\aa_{\alpha,\beta} \oplus \ll_{[x:y]}^{\alpha,\beta}) ;$$
\item soit il existe un unique $b \in N_{\alpha,\gamma}$ tel que $\ll = \Ad b(\aa_{\alpha,\gamma} \oplus \ll^{\alpha,\gamma})$ ;
\item soit il existe un unique $[x:y] \in \P^1(\R)$ et un unique $b \in N_{\beta,\gamma}$ tels que $$\ll = \Ad b(\aa_{\beta,\gamma} \oplus \ll_{[x:y]}^{\beta,\gamma}) ;$$
\item soit il existe un unique $[x:y:z] \in \P^2(\R)$ avec $x$ et $z$ non nuls tel que $\ll=\ll_{[x:y:z]}$ ;
\item soit il existe un unique $(y,t) \in \R^2$ tel que $\ll=\ll_{\alpha,y,t}$ ;
\item soit il existe un unique $(y,t) \in \R^2$ tel que $\ll=\ll_{\gamma,y,t}$ ;
\item soit il existe un unique $[x:y:z:t] \in \P^3(\R)$ tel que $\ll=\ll_{[x:y:z:t]}$ ;
\item soit il existe un unique $(x,y) \in \R^2$ tel que $\ll=\ll_{x,y}$.
\end{enumerate}
\epro

Ces sous-algèbres de Lie abéliennes sont toutes maximales, sauf $\ll_{[x:y:z:t]}$ qui est incluse dans $\ll_0$.

\bp Distinguons selon $p_\aa(\ll)$ et $p_{\alpha,\gamma}(\ll)$.
\ben
\item Si $p_\aa(\ll)=\aa$, considérons une base de $\ll$ constituée de trois éléments diagonalisables. Alors $\ll$ est un sous-espace de Cartan de $\gg$, il existe donc $b \in G$ tel que $\ll = \Ad b (\aa)$. Puisque $\bb$ est l'algèbre de Lie d'un sous-groupe de Borel contenant $\ll$, on peut supposer de plus que $\Ad b(\bb) = \bb$, c'est-à-dire que $b \in N_G(\bb)=B$. Puisque $MA$ normalise $A$, on peut supposer que $b \in N$. Et cet élément est unique car $N \cap N_G(A) = \{e\}$.
\item Si $p_\aa(\ll)$ est de dimension $2$, remarquons que si $X \in \bb$ est tel que $\alpha(p_\aa(X))$, $\beta(p_\aa(X))$ et $\gamma(p_\aa(X))$ soient tous les trois non nuls, alors le centralisateur de $X$ dans $\nn$ est trivial, donc par un argument de dimension $X$ n'appartient pas à $\ll$. Ainsi il existe une unique racine $\delta \in \Delta$ telle que $p_\aa(\ll)=\aa_\delta$. Alors la projection $\ll'$ de $\ll$ sur la sous-algèbre de Lie $\zz^{\Delta \bs \{\delta\}}$ est abélienne diagonalisable de dimension au moins $2$, donc il existe $b' \in N \cap Z^{\Delta \bs \{\delta\}}$ tel que $\Ad b'(\aa^{\Delta \bs \{\delta\}}) = \ll'$.
\bit
\item Si $\delta=\alpha$, il existe deux réels $x$ et $y$ tels que $(H_\beta + x U_{\alpha+\beta} + yU_{\alpha+\beta+\gamma},H_\gamma + yU_{\alpha+\beta+\gamma},U_\alpha)$ soit une base de $\Ad b'^{-1}(\ll)$. Posons $b''=\exp(xU_{\alpha+\beta}+yU_{\alpha+\beta+\gamma})$, alors l'élément $b=b'b'' \in N_\alpha$ est tel que $\ll = \Ad b(\aa_\alpha \oplus \nn^\alpha)$. Par ailleurs $b$ est unique, car le normalisateur de $\aa_\alpha \oplus \nn^\alpha$ dans $N$ est égal à $AN^\alpha$, et $AN^\alpha \cap N_\alpha = \{e\}$.
\item Si $\delta=\alpha$, il existe trois réels $x$, $y$ et $z$ tels que $(H_\alpha + x U_{\alpha+\beta} + zU_{\alpha+\beta+\gamma},H_\gamma + yU_{\beta+\gamma} + zU_{\alpha+\beta+\gamma},U_\beta)$ soit une base de $\Ad b'^{-1}(\ll)$. Définissons $b''=\exp(xU_{\alpha+\beta}+yU_{\beta+\gamma}+zU_{\alpha+\beta+\gamma}) \in N_\beta$, alors l'élément $b=b'b'' \in N_\beta$ est tel que $\ll = \Ad b(\aa_\beta \oplus \nn^\beta)$. Par ailleurs $b$ est unique, car le normalisateur de $\aa_\beta \oplus \nn^\beta$ dans $N$ est égal à $AN^\beta$, et $AN^\beta \cap N_\beta = \{e\}$.
\item Si $\delta=\alpha$, on montre comme dans le cas $\delta=\alpha$ qu'il existe un unique $b \in N_\gamma$ tel que $\ll = \Ad b(\aa_\gamma \oplus \nn^\gamma)$.
\eit
\item Si $p_\aa(\ll)$ est de dimension $1$, remarquons que si $X \in \bb$ est tel que deux des trois réels $\alpha(p_\aa(X))$, $\beta(p_\aa(X))$ et $\gamma(p_\aa(X))$ soient non nuls, alors le centralisateur de $X$ dans $\nn$ est de dimension $1$, donc $X$ n'appartient pas à $\ll$. Ainsi commençons par le cas où $p_\aa(\ll) = \aa_{\alpha,\beta}$, alors la projection $\ll'$ de $\ll$ sur la sous-algèbre de Lie $\zz^{\alpha,\beta} \simeq \sl_3(\R)$ est abélienne, donc de dimension au plus $2$. Soit $X \in \ll$ dont la projection sur $\zz^{\alpha,\beta}$ soit nulle, et dont la projection sur $\aa$ soit $H_\gamma$ : alors l'élément $b = \exp(X-H_\gamma) \in N_{\alpha,\beta}$ est tel que $X= \Ad b (H_\gamma)$. Or le centralisateur de $H_\gamma$ dans $\nn$ est $\nn^{\alpha,\beta} = \nn \cap \zz^{\alpha,\beta}$, donc la sous-algèbre $\ll'$ de $\zz^{\alpha,\beta} \simeq \sl_3(\R)$ est abélienne, de dimension $2$, incluse dans $\nn$ : d'après la proposition~\ref{pro:ssalg}, il existe un unique $[x:y] \in \P^1(\R)$ tel que $\ll' = \ll_{[x:y]}^{\alpha,\beta}$. Puisque l'élément $b$ normalise $\ll'$, on en déduit que $\ll = \Ad b(\aa_{\alpha,\beta} \oplus \ll_{[x:y]}^{\alpha,\beta})$. Par ailleurs $b$ est unique, car le normalisateur de $\aa_{\alpha,\beta} \oplus \ll_{[x:y]}^{\alpha,\beta}$ dans $N$ est égal à $N^{\alpha,\beta}$, et $N^{\alpha,\beta} \cap N_{\alpha,\beta} = \{e\}$.
\item Si $p_\aa(\ll) = \aa_{\alpha,\gamma}$, alors la projection $\ll'$ de $\ll$ sur la sous-algèbre de Lie $\zz^{\alpha,\gamma} \simeq \sl_2(\R) \oplus \sl_2(\R)$ est abélienne, donc de dimension au plus $2$. Soit $X \in \ll$ dont la projection sur $\zz^{\alpha,\gamma}$ soit nulle, et dont la projection sur $\aa$ soit $H_\beta$ : alors l'élément $b = \exp(X-H_\beta) \in N_{\alpha,\gamma}$ est tel que $X= \Ad b (H_\beta)$. Or le centralisateur de $H_\beta$ dans $\nn$ est $\nn^{\alpha,\gamma} = \nn \cap \zz^{\alpha,\gamma}$, donc la sous-algèbre $\ll'$ de $\zz^{\alpha,\gamma} \simeq \sl_2(\R) \oplus \sl_2(\R)$ est abélienne, de dimension $2$, incluse dans $\nn$, donc $\ll' = \ll^{\alpha,\gamma}$. Puisque l'élément $b$ normalise $\ll'$, on en déduit que $\ll = \Ad b(\aa_{\alpha,\gamma} \oplus \ll^{\alpha,\gamma})$. Par ailleurs $b$ est unique, car le normalisateur de $\aa_{\alpha,\gamma} \oplus \ll^{\alpha,\gamma}$ dans $N$ est égal à $N^{\alpha,\gamma}$, et $N^{\alpha,\gamma} \cap N_{\alpha,\gamma} = \{e\}$.
\item Si $p_\aa(\ll) = \aa_{\beta,\gamma}$, on montre comme dans le cas où $p_\aa(\ll) = \aa_{\alpha,\beta}$ qu'il existe un unique $[x:y] \in \P^1(\R)$ et un unique $b \in N_{\beta,\gamma}$ tels que $\ll = \Ad b(\aa_{\beta,\gamma} \oplus \ll_{[x:y]}^{\beta,\gamma})$.
\item Si $p_\aa(\ll)=\{0\}$ et $p_{\alpha,\gamma}(\ll)$ est de dimension $1$, soit $X \in \ll$ tel que $p_{\alpha,\gamma}(X) \neq 0$, et soient $x,y,z \in \R$ tels que $p_\alpha(X) = xU_\alpha$, $p_\beta(X) = yU_\beta$ et $p_\gamma(X) = zU_\gamma$. Supposons ici que $x$ et $z$ sont non nuls, alors le centralisateur de $X$ dans $\nn$ est $\ll_{[x:y:z]}$, donc $\ll=\ll_{[x:y:z]}$. Par ailleurs $[x:y:z]$ est unique.
\item Si $p_\aa(\ll)=\{0\}$ et $p_{\alpha,\gamma}(\ll)$ est de dimension $1$, soit $X \in \ll$ tel que $p_{\alpha,\gamma}(X) \neq 0$, et soient $x,y,z \in \R$ tels que $p_\alpha(X) = xU_\alpha$, $p_\beta(X) = yU_\beta$ et $p_\gamma(X) = zU_\gamma$. Supposons ici que $x=1$ et $z=0$, alors le centralisateur de $X$ dans $\nn$ est $\ll_{\alpha,y,t}$, où $t \in \R$ est tel que $p_{\beta+\gamma}(X) = tU_{\beta+\gamma}$. Ainsi $\ll=\ll_{\alpha,y,t}$, et $(y,t) \in \R^2$ est unique.
\item Si $p_\aa(\ll)=\{0\}$ et $p_{\alpha,\gamma}(\ll)$ est de dimension $1$, soit $X \in \ll$ tel que $p_{\alpha,\gamma}(X) \neq 0$, et soient $x,y,z \in \R$ tels que $p_\alpha(X) = xU_\alpha$, $p_\beta(X) = yU_\beta$ et $p_\gamma(X) = zU_\gamma$. Supposons ici que $x=0$ et $z=1$, alors le centralisateur de $X$ dans $\nn$ est $\ll_{\gamma,y,t}$, où $t \in \R$ est tel que $p_{\alpha+\beta}(X) = tU_{\alpha+\beta}$. Ainsi $\ll=\ll_{\gamma,y,t}$, et $(y,t) \in \R^2$ est unique.
\item Si $p_\aa(\ll)=\{0\}$ et $p_{\alpha,\gamma}(\ll)$ est de dimension $0$, alors $\ll$ est incluse dans la sous-algèbre de Lie abélienne $\ll_0$ de dimension $4$, donc il existe un unique $[x:y:z:t] \in \P^3(\R)$ tel que $\ll=\ll_{[x:y:z:t]}$.
\item Si $p_\aa(\ll)=\{0\}$ et $p_{\alpha,\gamma}(\ll)$ est de dimension $2$, considérons $X,Y \in \ll$ tels que $p_{\alpha,\gamma}(X)=U_\alpha$ et $p_{\alpha,\gamma}(Y)=U_\gamma$. Soit $(x,y) \in \R^2$ tels que $p_{\beta+\gamma}(X)=xU_{\beta+\gamma}$ et $p_{\alpha+\beta}(Y)=yU_{\alpha+\beta}$. Alors le centralisateur de $X$ et $Y$ dans $\nn$ est $\ll_{x,y}$, de dimension $3$, donc $\ll=\ll_{x,y}$. Et l'élément $(x,y) \in \R^2$ est unique.
\een
\epp

Nous disposons alors d'un théorème analogue au théorème~\ref{thm:sgabeliens_sl3} pour $\SL_4(\R)$.

\bthm \label{thm:sgabeliens_sl4} La compactification de Chabauty $\ov{\Cartan(\SL_4(\R))}^\Ch$ coïncide avec l'espace ${\cal A}(\SL_4(\R))$. \ethm

\bp Il suffit de montrer que toutes les sous-algèbres de Lie décrites dans la proposition~\ref{pro:ssalg_sl4} sont limites de sous-espaces de Cartan. Soit $\ll$ l'une de ces sous-algèbres de Lie.
\ben
\item Si $\ll = \aa$, il n'y a rien à démontrer.
\item S'il existe une racine primitive $\delta \in \Delta$ telle que $\ll = \aa_\delta \oplus \nn^\delta$, alors la suite de sous-espaces de Cartan $(\Ad \exp(nU_\delta) (\aa))_{n \in \N}$ converge vers $\ll$.
\item S'il existe $[x:y] \in \P^1(\R)$ tel que $\ll = \aa_{\alpha,\beta} \oplus \ll_{[x:y]}^{\alpha,\beta}$, alors d'après le théorème~\ref{thm:sgabeliens_sl3} appliqué au sous-groupe $Z^{\alpha,\beta} \simeq \SL_3(\R)$, la sous-algèbre de Lie $\ll_{[x:y]}^{\alpha,\beta}$ de $\zz^{\alpha,\beta}$ (dont le sous-espace de Cartan standard est $\aa^{\alpha,\beta}$) est la limite d'une suite $(\Ad g_n(\aa^{\alpha,\beta}))_{n \in \N}$, où $g_n \in Z^{\alpha,\beta}$ pour tout $n \in \N$. Alors la suite $(\Ad g_n(\aa))_{n \in \N}$ converge vers $\ll$.
\item Si $\ll = \aa_{\alpha,\gamma} \oplus \ll_{\alpha,\gamma}$, alors la suite $(\Ad \exp(nU_\alpha+nU_\gamma) (\aa))_{n \in \N}$ converge vers $\ll$.
\item S'il existe $[x:y] \in \P^1(\R)$ tel que $\ll = \aa_{\beta,\gamma} \oplus \ll_{[x:y]}^{\beta,\gamma}$, alors comme pour le cas $\aa_{\alpha,\beta} \oplus \ll_{[x:y]}^{\alpha,\beta}$, on montre que $\ll$ est limite de sous-espaces de Cartan.
\item S'il existe $[x:y:z] \in \P^2(\R)$ avec $x$ et $z$ non nuls tel que $\ll=\ll_{[x:y:z]}$, distinguons deux cas.
\bit
\item Si $y$ est non nul, alors quitte à conjuguer $\ll$ par un élément de $A$ on peut supposer que $[x:y:z]=[1:1:1]$. Toute valeur d'adhérence de la suite de sous-espaces de Cartan $(\Ad \exp(nU_\alpha+nU_\beta+nU_\gamma) (\aa))_{n \in \N}$ contient le vecteur
$$\liml_{n \ra +\infty} \Ad \exp(nU_\alpha+nU_\beta+nU_\gamma)\left(-\frac{1}{n}H_\alpha-\frac{1}{n}H_\beta-\frac{1}{n}H_\gamma\right) = U_\alpha + U_\beta + U_\gamma,$$
et par ailleurs toute valeur d'adhérence de cette suite est incluse dans $\nn$, donc d'après la proposition~\ref{pro:ssalg_sl4} nous en déduisons que cette suite de sous-espaces de Cartan converge vers $\ll$.
\item Si $y=0$, alors la suite de sous-algèbres de Lie $\left(\ll_{[x:\frac{1}{n}:z]}\right)_{n \in \N \bs \{0\}}$ converge vers $\ll=\ll_{[x:y:z]}$. Ainsi d'après le point précédent $\ll$ est limite de sous-espaces de Cartan.
\eit
\item S'il existe $(y,t) \in \R^2$ tel que $\ll=\ll_{\alpha,y,t}$, distinguons deux cas.
\bit
\item Si $y$ et $t$ sont non nuls, alors quitte à conjuguer par $A$ on peut supposer que $y=t=1$. Pour tout $n \in \N$, posons $b_n=\exp\left(nU_\alpha+nU_\beta+nU_{\beta+\gamma}\right) \in N$. Alors toute valeur d'adhérence de la suite $(\Ad b_n(\aa))_{n \in \N}$ contient le vecteur
$$\liml_{n \ra +\infty} \Ad b_n \left(-\frac{1}{n}H_\alpha-\frac{1}{n}H_\beta \right) = U_\alpha + U_\beta + U_{\beta+\gamma},$$
et par ailleurs toute valeur d'adhérence de cette suite est incluse dans $\nn$, donc d'après la proposition~\ref{pro:ssalg_sl4} nous en déduisons que cette suite de sous-espaces de Cartan converge vers $\ll$.
\item Si $y$ ou $t$ est nul, alors considérons deux suites réelles $(y_n)_{n \in \N}$ et $(t_n)_{n \in \N}$, dont tous les termes sont non nuls, convergeant vers $y$ et $t$ respectivement. Alors la suite de sous-algèbres de Lie $(\ll_{\alpha,y_n,t_n})_{n \in \N}$ converge vers $\ll=\ll_{\alpha,y,t}$. Ainsi d'après le point précédent $\ll$ est limite de sous-espaces de Cartan.
\eit
\item S'il existe $(y,t) \in \R^2$ tel que $\ll=\ll_{\gamma,y,t}$, alors comme dans le cas précédent on montre que $\ll$ est limite de sous-espaces de Cartan.
\item S'il existe $[x:y:z:t] \in \P^3(\R)$ tel que $\ll=\ll_{[x:y:z:t]}$, alors considérons l'action par conjugaison du sous-groupe
$$ N_G(\ll_0) = \left\{ \left( \begin{array}{cccc} *&*&*&* \\ *&*&*&* \\ 0&0&*&* \\ 0&0&*&* \end{array} \right) \right\}.$$
Si le rang de la matrice $\left( \begin{array}{cc} x&y \\ z&t \end{array} \right)$ est $2$ on peut supposer que $[x:y:z:t] = [1:0:0:1]$, et si le rang est $1$ on peut supposer que $[x:y:z:t] = [0:0:1:0]$.
\bit
\item Si $[x:y:z:t] = [1:0:0:1]$, alors remarquons que la suite $(\ll=\ll_{[1:n:1]})_{n \in \N}$ converge vers $\ll$, or d'après le point précédent nous savons que, pour tout $n \in \N$, la sous-algèbre de Lie $\ll_{[1:n:1]}$ est limite de sous-espaces de Cartan, donc c'est également le cas de $\ll$.
\item Si $[x:y:z:t] = [0:0:1:0]$, alors remarquons que la suite $(\ll=\ll_{[1:0:n:1]})_{n \in \N}$ converge vers $\ll$, or d'après le point précédent nous savons que, pour tout $n \in \N$, la sous-algèbre de Lie $\ll_{[1:0:n:1]}$ est limite de sous-espaces de Cartan, donc c'est également le cas de $\ll$.
\eit
\item S'il existe $(x,y) \in \R^2$ tel que $\ll=\ll_{x,y}$, distinguons deux cas.
\bit
\item Si $x$ et $y$ sont non nuls, alors quitte à conjuguer par $A$ on peut supposer que $x=y=-1$. Soit $b_n=\exp\left(nU_\alpha+nU_\gamma+nU_{\alpha+\beta}+nU_{\beta+\gamma}+n^4U_{\alpha+\beta+\gamma}\right) \in N$, pour tout $n \in \N$. Alors toute valeur d'adhérence de la suite $(\Ad b_n(\aa))_{n \in \N}$ contient les vecteurs
\beq \liml_{n \ra +\infty} \Ad b_n \left(\frac{1}{n}H_\alpha-\frac{1}{n}H_\beta+\frac{1}{n^3}H_\gamma \right) &=& U_\alpha - U_{\beta+\gamma} \\
\mbox{et } \liml_{n \ra +\infty} \Ad b_n \left(-\frac{1}{n^3}H_\alpha-\frac{1}{n}H_\beta+\frac{1}{n}H_\gamma \right) &=& U_\gamma - U_{\alpha+\beta},\eeq
et par ailleurs toute valeur d'adhérence de cette suite est incluse dans $\nn$, donc d'après la proposition~\ref{pro:ssalg_sl4} nous en déduisons que cette suite de sous-espaces de Cartan converge vers $\ll$.
\item Si $x$ ou $y$ est nul, alors considérons deux suites réelles $(x_n)_{n \in \N}$ et $(y_n)_{n \in \N}$, dont tous les termes sont non nuls, convergeant vers $x$ et $y$ respectivement. Alors la suite de sous-algèbres de Lie $(\ll_{x_n,y_n})_{n \in \N}$ converge vers $\ll=\ll_{x,y}$. Ainsi d'après le point précédent $\ll$ est limite de sous-espaces de Cartan.
\eit
\een
Nous avons donc montré que toutes les sous-algèbres de Lie de $\gg$ abéliennes de dimension $3$ incluses dans $\bb$ sont limites de sous-espaces de Cartan. Puisque $K$ agit transitivement sur les sous-algèbres de Borel de $\sl_4(\R)$, ceci montre que $\ov{\Cartan(\sl_4(\R))}^\Ch = {\cal A}(\sl_4(\R))$.
\ep

\bibliographystyle{smfalpha_perso}
\bibliography{biblio}

\providecommand{\bysame}{\leavevmode ---\ }
\providecommand{\og}{``}
\providecommand{\fg}{''}
\providecommand{\smfandname}{\&}
\providecommand{\smfedsname}{\'eds.}
\providecommand{\smfedname}{\'ed.}
\providecommand{\smfmastersthesisname}{M\'emoire}
\providecommand{\smfphdthesisname}{Th\`ese}
\begin{thebibliography}{LBG11b}

\bibitem[All99]{allcock_cayley}
{\scshape D.~Allcock} -- {\og {Reflections Groups on the Octave Hyperbolic
  Plane}\fg}, \emph{{J. Algebra}} \textbf{213} (1999), p.~467--498.

\bibitem[BCR87]{bochnak_coste_roy}
{\scshape J.~Bochnak, M.~Coste {\normalfont \smfandname} M.-F. Roy} --
  \emph{{Géométrie algébrique réelle\!}} , {Ergebnisse der Mathematik und ihrer
  Grenzgebiete (3)}, Springer-Verlag, 1987.

\bibitem[CG74]{chen_greenberg}
{\scshape S.~Chen {\normalfont \smfandname} L.~Greenberg} -- \emph{{Hyperbolic
  spaces}} , \dans \og Contributions to Analysis \fg, Academic Press, New York,
  1974.

\bibitem[Cha50]{chabauty}
{\scshape C.~Chabauty} -- {\og {Limite d'ensembles et g\'eom\'etrie des
  nombres}\fg}, \emph{{Bull. Soc. Math. France}} \textbf{78} (1950),
  p.~143--151.

\bibitem[CL09]{caprace_lecureux}
{\scshape P.-E. Caprace {\normalfont \smfandname} J.~Lécureux} -- {\og
  Combinatorial and group-theoretic compactifications of buildings\fg},
  (2009), {à paraître dans Ann. Inst. Fourier}.

\bibitem[dlH08]{harpe_chabauty}
{\scshape P.~de~la Harpe} -- {\og {Spaces of closed subgroups of locally
  compact groups}\fg}, {arXiv:0807.2030v2}, 2008.

\bibitem[FM94]{fulton_macpherson}
{\scshape W.~Fulton {\normalfont \smfandname} R.~MacPherson} -- {\og {A
  Compactification of Configuration Spaces}\fg}, \emph{{Annals of Math.}}
  \textbf{139} (1994), p.~183--225.

\bibitem[GJT98]{guivarch}
{\scshape Y.~Guivarc'h, L.~Ji {\normalfont \smfandname} J.~C. Taylor} --
  \emph{{Compactifications of symmetric spaces\!}} , {Progr.~Math. {\bf 156}},
  {Birkh\"auser}, 1998.

\bibitem[GR06]{guivarch_remy}
{\scshape Y.~Guivarc'h {\normalfont \smfandname} B.~Rémy} -- {\og
  {Group-theoretic compactification of Bruhat-Tits buildings}\fg}, \emph{Ann.
  Sci. École Norm. Sup.} \textbf{39} (2006), p.~871--920.

\bibitem[Hae10]{haettel_m2}
{\scshape T.~Haettel} -- {\og {Compactification de Chabauty des espaces
  symétriques de type non compact}\fg}, \emph{{J. Lie Theory}} \textbf{20}
  (2010), p.~437--468.

\bibitem[Hum75]{humphreys}
{\scshape J.~E. Humphreys} -- \emph{{Linear algebraic groups\!}} , {Grad. Texts
  Math. {\bf 21}}, {Springer Verlag}, 1975.

\bibitem[IM05a]{iliev_manivel_severi}
{\scshape A.~Iliev {\normalfont \smfandname} L.~Manivel} -- {\og {Severi
  varieties and their varieties of reductions}\fg}, \emph{{J. de Crelle}}
  \textbf{585} (2005), p.~93--139.

\bibitem[IM05b]{iliev_manivel}
\bysame , \emph{{Varieties of reduction for gln}} , \dans \og Projective
  Varieties with Unexpected Properties \fg, Walter de Gruyter, 2005.

\bibitem[LBG11a]{lebarbier_ex}
{\scshape M.~Le~Barbier~Grünewald} -- {\og {Examples of varieties of reductions
  of small rank}\fg}, {http://www.uni-bonn.de/~mlbg/public/michi-redex.pdf},
  {2011}.

\bibitem[LBG11b]{lebarbier}
\bysame , {\og {The variety of reductions for a reductive symmetric pair}\fg},
  \emph{{Transformation Groups}} \textbf{16} (2011), p.~1--26.

\bibitem[Mar97]{margulis_conjecture}
{\scshape G.~Margulis} -- {\og {Oppenheim conjecture}\fg}, \emph{Fields
  Medallists' lectures, \dans World Sci. Ser. 20th Century Math., v.5, World
  Sci. Publ.} (1997), p.~272--327.

\bibitem[Mau10]{maucourant}
{\scshape F.~Maucourant} -- {\og A non-homogeneous orbit closure of a diagonal
  subgroup\fg}, \emph{Ann. of Math.} \textbf{171} (2010), p.~557--570.

\bibitem[Moo79]{moore_amenable}
{\scshape C.~C. Moore} -- {\og Amenable subgroups of semi-simple groups and
  proximal flows\fg}, \emph{Isr. J. Math.} \textbf{34} (1979), p.~121--138.

\bibitem[Par07]{parker_hyperbolic}
{\scshape J.~R. Parker} -- {\og {Hyperbolic Spaces}\fg}, \emph{{The Jyväskylä
  Notes}} (2007).

\bibitem[Sat60]{satake}
{\scshape I.~Satake} -- {\og {On representations and compactifications of
  symmetric Riemannian spaces}\fg}, \emph{Ann. of Math.} \textbf{71} (1960),
  p.~77--110.

\bibitem[WZ84]{winternitz_zassenhaus}
{\scshape P.~Winternitz {\normalfont \smfandname} H.~Zassenhaus} -- \emph{{The
  structure of maximal abelian subalgebras of classical Lie and Jordan
  algebras\!}} , \dans \og XIIIth international colloquium on group theoretical
  methods in physics (College Park, Md., 1984) \fg, World Sci. Publishing,
  1984.

\end{thebibliography}

\sign

\end{document}